\documentclass[preprint,10pt]{elsarticle}
\usepackage[margin=1in]{geometry}

\usepackage{graphicx} 
\usepackage{bm}
\usepackage{amsfonts,amsmath}
\usepackage[rightcaption]{sidecap}
\usepackage{amsthm}
\usepackage{thmtools,thm-restate}

\usepackage{appendix}

\usepackage{mathtools}
\usepackage{xcolor}
\usepackage{natbib}
\usepackage{todonotes}
\usepackage[normalem]{ulem}

\usepackage{algorithm2e}
\usepackage{algpseudocode}
\algrenewcommand\algorithmicrequire{\textbf{Input:}}
\algrenewcommand\algorithmicensure{\textbf{Output:}}

\newtheoremstyle{DefStyle}
  {10pt} 
  {5pt} 
  {} 
  {} 
  {\bfseries} 
  {.} 
  {10pt} 
  {} 
\theoremstyle{DefStyle}

\newtheorem{theorem}[algocf]{Theorem} 

\newtheorem{definition}[theorem]{Definition}
\newtheorem{proposition}[theorem]{Proposition}

\let\oldproofname=\proofname
\renewcommand{\proofname}{\rm\bf{\oldproofname}}

\newcommand{\Rbb}{\mathbb{R}}

\newcommand{\Dc}{\mathcal{D}}
\newcommand{\Ic}{\mathcal{I}}

\newcommand{\Vc}{\mathcal{V}}

\newcommand{\Tsf}{\mathsf{T}}

\newcommand{\Absf}{\bm{\mathsf{A}}}
\newcommand{\Bbsf}{\bm{\mathsf{B}}}
\newcommand{\Cbsf}{\bm{\mathsf{C}}}
\newcommand{\Dbsf}{\bm{\mathsf{D}}}
\newcommand{\Ebsf}{\bm{\mathsf{E}}}
\newcommand{\Ibsf}{\bm{\mathsf{I}}}
\newcommand{\Lbsf}{\bm{\mathsf{L}}}
\newcommand{\Kbsf}{\bm{\mathsf{K}}}
\newcommand{\Mbsf}{\bm{\mathsf{M}}}
\newcommand{\Nbsf}{\bm{\mathsf{N}}}
\newcommand{\Pbsf}{\bm{\mathsf{P}}}
\newcommand{\Qbsf}{\bm{\mathsf{Q}}}
\newcommand{\Ubsf}{\bm{\mathsf{U}}}
\newcommand{\Rbsf}{\bm{\mathsf{R}}}
\newcommand{\absf}{\bm{\mathsf{a}}}
\newcommand{\Sbsf}{\bm{\mathsf{S}}}

\newcommand{\Wbsf}{\bm{\mathsf{W}}}
\newcommand{\Xbsf}{\bm{\mathsf{X}}}
\newcommand{\Ybsf}{\bm{\mathsf{Y}}}

\newcommand{\Itbsf}{\widetilde{\Ibsf}}
\newcommand{\Ytbsf}{\widetilde{\Ybsf}}
\newcommand{\ytbsf}{\widetilde{\ybsf}}

\newcommand{\cbsf}{\bm{\mathsf{c}}}
\newcommand{\ebsf}{\bm{\mathsf{e}}}

\newcommand{\nbsf}{\bm{\mathsf{n}}}

\newcommand{\qbsf}{\bm{\mathsf{q}}}
\newcommand{\ubsf}{\bm{\mathsf{u}}}
\newcommand{\vbsf}{\bm{\mathsf{v}}}
\newcommand{\wbsf}{\bm{\mathsf{w}}}
\newcommand{\xbsf}{\bm{\mathsf{x}}}
\newcommand{\ybsf}{\bm{\mathsf{y}}}

\newcommand{\ybsfs}{\ybsf^{\star}}

\newcommand{\Lambdabsf}{\bm{\mathsf{\Lambda}}}

\DeclareMathOperator*{\argmin}{argmin}
\newcommand{\DIAG}{{\sf diag}}
\newcommand{\TRACE}{{\sf Tr}}

\newcommand{\SUCHTHAT}{~:~}

\newcommand{\IdentityMatrix}{\Ibsf}
\newcommand{\ZeroMatrix}{\bm{\mathsf{0}}}

\newcommand{\ag}[1]{\left[#1 \right ]}
\renewcommand{\u}[1]{\boldsymbol{#1}}

\newcommand{\usf}[1]{\u{\mathsf #1}}
\newcommand{\busf}[1]{\overline{\u{\mathsf #1}}}

\newcommand{\pr}[1]{\left( #1 \right)}

\newcommand{\norm}[1]{\left\lVert #1\right\rVert}
\newcommand{\fprod}[1]{\langle #1 \rangle_{F}}

\newcommand{\pder}[2]{\frac{\partial #1}{\partial #2}}

\newcommand{\lsc}[1]{{}^{#1}\negthickspace\thinspace}

\newcommand{\Mvspace}[1]{\mathcal{M}_{#1}}
\newcommand{\Mprime}[1]{\widetilde{\mathcal{M}}_{#1}}
\newcommand{\Mpprod}[1]{\langle #1 \rangle_{\sim}}
\newcommand{\Kvspace}[1]{\mathcal{K}_{#1}}
\newcommand{\Svspace}[1]{\mathcal{S}_{#1}}

\newcommand{\Oset}[1]{\mathcal{O}_{#1}}
\newcommand{\Oprime}[1]{\widetilde{\mathcal{O}}_{#1}}
\newcommand{\Ktwoset}[1]{\mathcal{K}^{\rm sq}_{#1}}

\usepackage{amsmath}

\usepackage{tikz}

\usepackage{mathabx}

\usepackage{ifthen}
\newcommand{\ShowComments}{false}
\newcommand{\forus}[1]{\ifthenelse{\equal{\ShowComments}{true}}{{\footnotesize \color{blue!90!black}~[[#1]]}}{}}

\makeatletter
\newcommand\footnoteref[1]{\protected@xdef\@thefnmark{\ref{#1}}\@footnotemark}
\makeatother

\usepackage{longtable}
\usepackage{tablefootnote}
\usepackage{tabularray}
\UseTblrLibrary{booktabs}

\begin{document}

\title{Approximating a matrix as the square of a skew-symmetric matrix, with application to estimating angular velocity from acceleration data}

\author[1]{Yang Wan}
\author[1,2]{Benjamin E.~Grossman-Ponemon}
\author[1]{Haneesh Kesari\corref{cor1}}
\ead{haneesh\_kesari@brown.edu}
\address[1]{School of Engineering, Brown University, Providence, RI 02912, USA}
\address[2]{New address: Department of Physics and Engineering, John Carroll University, University Heights, OH 44118, USA}
\cortext[cor1]{Corresponding author}

\begin{abstract}
  In this paper we study the problem of finding the best approximation of a real square matrix by a matrix that can be represented as the square of a real, skew-symmetric matrix. %
  This problem is important in the design of robust numerical algorithms aimed at estimating rigid body kinematics from multiple accelerometer measurements.
  We give a constructive proof for the existence of a best approximant in the Frobenius norm.
  We demonstrate the construction with some small examples, and we showcase the practical importance of this work to the problem of determining the angular velocity of a rotating rigid body from its acceleration measurements.
\end{abstract}

\maketitle

\section{Introduction}

Approximating a given matrix by one with special properties appears in a variety of contexts in applied mathematics and engineering. For example, the problems of approximating real, square matrices in the set of symmetric, positive, semidefinite matrices~\cite{Halmos1972,hayden1988approximation,higham1988computing} and in the orthonormal group~\cite{Higham1986} see applications in computational methods for machine learning and control theory.

Here, we concern ourselves with the following matrix approximation problem. Given a real, square matrix $\Absf$, can we find a matrix $\Ubsf = \Kbsf^2$---where $\Kbsf$ is a real, skew-symmetric matrix---that best approximates $\Absf$ among all such matrices? More precisely, if $\Mvspace{n}$ is the Euclidean space of real, $n\times n$ matrices equipped with the Frobenius inner product, and $\Kvspace{n} \subset \Mvspace{n}$ is the subspace of real, skew-symmetric, $n\times n$ matrices, then we solve
\begin{equation}
    \argmin_{\Ubsf \in \Ktwoset{n}} \|\Ubsf - \Absf\|_F,
    \tag{$\mathbb{P}$}
    \label{p:OptimizationProblem}
\end{equation}
where
\begin{equation}
\Ktwoset{n} = \left\{ \Kbsf^2 \SUCHTHAT \Kbsf \in \Kvspace{n} \right\}
\label{def:ProblemDomain}
\end{equation}
and $\| \cdot \|_F$ denotes the norm induced by the Frobenius inner product. As the above minimization depends on the matrix $\Absf \in \Mvspace{n}$, we will refer to the problem as Problem \ref{p:OptimizationProblem}($\Absf$).

The main contribution of this work is an explicit construction for a minimizer of Problem \ref{p:OptimizationProblem}($\Absf$). We remark that this may not be the sole minimizer of Problem \ref{p:OptimizationProblem}($\Absf$). We provide the construction in two formats: first as Definition \ref{def:S3Approx}, and later as Algorithm \ref{algo:ustar}. Because our minimizer is provided explicitly, we do not need to use optimization techniques such as Newton's Method to solve Problem \ref{p:OptimizationProblem}($\Absf$) which may only find local minima.

We are not the first to consider this optimization problem. In the case where $n = 3$, Problem \ref{p:OptimizationProblem}($\Absf$) has been attempted in the literature~\cite{Cardou2011}. However, the construction from this previous work is obtained by investigating the stationary points of Problem \ref{p:OptimizationProblem}($\Absf$)\footnote{To be precise, there exists a map from $\Rbb^3$ to $\Kvspace{3}$ (and hence $\Ktwoset{3}$). In Ref.~\cite{Cardou2011}, Problem \ref{p:OptimizationProblem}($\Absf$) is reformulated as a minimization over $\Rbb^3$. It is for this problem that the stationary points are computed.}. The construction is not proved to produce a global minimizer; moreover, as presented in~\cite[Section III.B]{Cardou2011}, the construction does not work for certain input matrices $\Absf$.\footnote{These cases are accounted for in the computer algorithm presented in~\cite[Appendix B]{Cardou2011}, though it is unclear how these changes relate to the minimizers of the original optimization problem.}

To conclude this section, we elaborate on the immediate practical application of our study below. Afterwards, this work is structured as follows. A formal statement of the mathematical problem is given in \S\ref{sec:ProblemStatementAndSolution}, wherein our solution to Problem \ref{p:OptimizationProblem}($\Absf$) is posed as Theorem~\ref{thm:main}.  In preparation for the proof of Theorem~\ref{thm:main},  we recall a few identities in \S\ref{sub:SupportingResults}. The proof of Theorem~\ref{thm:main} then follows in \S\ref{sub:ProofOfMainTheorem}. We return to the practical application of our study in \S\ref{sec:exaapp}, illustrating Algorithm \ref{algo:ustar} with some examples in \S\ref{subsec:examples} and showcasing an experimental demonstration in \S\ref{sec:expdemo}. Concluding remarks are provided in \S\ref{sec:ConcludingRemarks}.

\subsection{Motivating application}\label{sec:MotivationApplication}
Angular velocity is a vector-valued measure of a rigid body's rate of rotation. Its experimental measurement is critical in a number of fields, such as aeronautics \cite{Orlik1983}, astronautics \cite{a15020029,MA201489,LIANG20111019,GUGLIERI2014395}, robotics \cite{BFb0030799,Barshan1995}, and, most recently, biomechanics \cite{takhounts2013development,Wan2023brain,Carlsen2021}. The head's angular velocity during a traumatic event, such as a blunt impact or a fall, is considered as a key parameter for assessing that event's  risk of leading to mild traumatic brain injury (mTBI) \cite{Wan2023brain,Carlsen2021}.

The $\sqrt{\text{AO}}$-algorithm~\cite{wan2022determining} is a method that estimates rigid body kinematics (including angular velocity and angular acceleration) using only measurements from accelerometers, rather than from gyroscopes.
This algorithm falls into the class of what are termed \textit{gyroscope-free algorithms} \cite{parsa2003dynamics,Cardou2008,Cardou2011,parsa2005design}.
Accelerometers typically have larger bandwidths than gyroscopes when they are of similar size and weight (Table.~\ref{tab:memssensor}), and so gyroscope-free algorithms allow for measurement of angular velocity at much higher frequencies since they rely only on accelerometer data.

We briefly outline the template of the $\sqrt{\rm AO}$-algorithm below.
The algorithm's inputs are measurements of the acceleration vectors at four of the rigid body's points at a discrete sequence of time instances. (See last paragraph of \S\ref{sec:accemeasurement} for elaboration.)
The acceleration vectors of the rigid body's points are usually measured in different bases (\textit{accelerometer bases}\footnote{e.g., the orthonormal vector sets $\pr{\lsc{\ell}\u{e}_i\ag{\tau}}_{i\in (1,2,3)}$, $\ell=1,2,3,4$, shown in Fig.~\ref{fig:testsetup}(b.ii)}).
The accelerometer bases are attached to the rigid body and move with it (cf. Figs.~\ref{fig:testsetup}(b.i)--(b.ii)). A plain version of the $\sqrt{\text{AO}}$-algorithm consists of the following two steps.
~\textit{(i)} At each time instance $\tau_i$, $i = 1,2,\ldots$, the acceleration data are synthesized to produce the symmetric matrix $\Bbsf\ag{\tau_i}$ (see Ref.~\cite{wan2022determining} for details).
The computation for $\Bbsf\ag{\tau_i}$ also requires information about the relative distances between the rigid body's points, and the relative orientations of the accelerometer bases.
~\textit{(ii)} The angular velocity matrix of the rigid body at the time instance $\tau_i$, $\overline{\Wbsf}\ag{\tau_i}\in \Kvspace{n}$, is then computed from the roots\footnote{
 When multiple roots exist, the physically meaningful root is taken to be the one that is closest to the physically meaningful root from the previous time instance. Thus, in effect, the rigid body's initial angular velocity dictates the choice of the physically meaningful root at each time instance.} of the function $r\ag{\Bbsf\ag{\tau_i}}\ag{\cdot}$, where $r\ag{\Bbsf}\ag{\cdot} : \Kvspace{n} \rightarrow \Svspace{n}$,
\begin{equation}
r\ag{\Bbsf}\ag{\Wbsf}=\Wbsf\,\Wbsf-\Bbsf,
\label{eq:symmp}
\end{equation}
and $\Bbsf\in  \Svspace{n}$,  the subspace of $\Mvspace{n}$ consisting of all real, symmetric, $n \times n$ matrices.

In theory, applying the plain $\sqrt{\text{AO}}$-algorithm to compute angular velocities is straightforward; in practice, however, this proves challenging.

The set $\Ktwoset{n}$ defined in \eqref{def:ProblemDomain} is a proper subset of $\Svspace{n}$, though not possessing the same vector space structure.
It can be shown that $r[\Bbsf]\ag{\cdot}$ must have a root when $\Bbsf \in \mathcal{K}_n^{\rm sq}$ and have no roots when $\Bbsf\in \mathcal{S}_n\setminus \mathcal{K}_n^{\rm sq}$. For each $\tau_i$, the matrix $\Bbsf[\tau_i]$ always belongs to $\Svspace{n}$.
When $\Bbsf[\tau_i]$ comes from perfect rigid body motion, then $\Bbsf[\tau_i] \in \Ktwoset{n}$.
Therefore, in theory, $\overline{\Wbsf}\ag{\tau_i}$ can always be computed from the roots of $r\ag{\Bbsf\ag{\tau_i}}\ag{\cdot}$ (cf.~\cite[\S 4.2]{wan2022determining}), i.e., by applying  step \textit{(ii)} of the plain $\sqrt{\text{AO}}$-algorithm.

Due to noise in real acceleration data and lack of precision in position and orientation measurements, applying step \textit{(i)} of the plain $\sqrt{\text{AO}}$-algorithm only yields an approximation for $\Bbsf\ag{\tau_i}$. Let us call this approximation $\widetilde{\Bbsf}\ag{\tau_i}$. Then what we can hope to get from $\widetilde{\Bbsf}\ag{\tau_i}$ is only an approximation for the actual angular velocity $\overline{\Wbsf}\ag{\tau_i}$. This is not an issue in itself, since most experimental procedures also only yield approximations for the physical quantity of interest. The issue lies in that $\widetilde{\Bbsf}\ag{\tau_i}$ rarely, if ever, belongs to $\Ktwoset{n}$, and so $r\ag{\widetilde{\Bbsf}\ag{\tau_i}}\ag{\cdot}$ has no roots. Consequently, step \textit{(ii)} of the plain $\sqrt{\text{AO}}$-algorithm does not work.

An alternate strategy for computing $\widetilde{\overline{\Wbsf}}\ag{\tau_i}$ from $\widetilde{\Bbsf}\ag{\tau_i}$ is based on finding best approximants as follows.  By construction, $\widetilde{\Bbsf}\ag{\tau_i}$ always belongs to $\Svspace{n}$.
Find a best approximant for $\widetilde{\Bbsf}\ag{\tau_i}$ in $\Ktwoset{n}$. Let this matrix be $\widehat{\Bbsf}\ag{\tau_i}$.  Compute the roots of $r\ag{\widehat{\Bbsf}\ag{\tau_i}}\ag{\cdot}$. This is now straightforward, since  $\widehat{\Bbsf}\ag{\tau_i}\in \Ktwoset{n}$.  Finally, compute $\widetilde{\overline{\Wbsf}}\ag{\tau_i}$ from the roots of $r\ag{\widehat{\Bbsf}\ag{\tau_i}}\ag{\cdot}$, in the same manner as $\overline{\Wbsf}\ag{\tau_i}$ is computed from the roots of $r\ag{\Bbsf\ag{\tau_i}}\ag{\cdot}$. This alternate strategy is the $\sqrt{\text{AO}}$-algorithm. Note that it is the same as the plain $\sqrt{\text{AO}}$-algorithm except for the additional approximation step, i.e., approximating $\widetilde{\Bbsf}\ag{\tau_i}$ with $\widehat{\Bbsf}\ag{\tau_i}$. Thus, success of the $\sqrt{\text{AO}}$-algorithm is dependent on this approximation step, namely approximating an arbitrary, real, symmetric matrix by one belonging to $\Ktwoset{n}$.

In the sequel, we outline the problem we study, which is a more general version of the problem needed for the $\sqrt{\text{AO}}$-algorithm. In this more general problem, we consider determining a best approximant in $\Ktwoset{n}$ of an arbitrary real, $n \times n$ matrix $\Absf$ rather than a real, symmetric  $n \times n$ matrix $\Bbsf$. Our analysis is general with respect to the dimension $n$, whereas for the $\sqrt{\text{AO}}$-algorithm it suffices to consider only the cases $n=2$ and $n=3$.

\section{Preliminaries and main result}
\label{sec:ProblemStatementAndSolution}


\subsection{Definitions and notation}
\label{sec:notation}

To make the statement of \ref{p:OptimizationProblem} more precise we first introduce some definitions and notation.

\paragraph{Matrix components:} 
Say the matrix $\Absf=\pr{\pr{A_{ij}}_{j\in \Ic}}_{i\in \Ic}$, where
$A_{ij}\in \mathbb{R}$, and $\Ic=(1,\ldots, n)$.
We sometimes denote the $p\mbox{-}q^{\rm th}$ component of $\Absf$, i.e.~$A_{pq}$, as $\Absf_{\cdot p\cdot q}$. 

\paragraph{Diagonal matrices:} 
An important subset of $\Svspace{n}$ is the set of diagonal matrices. 
 We will compactly write these matrices using the operator $\DIAG_n\ag{\cdot,\cdot,\ldots, \cdot} : \Rbb^n \rightarrow \Svspace{n}$,
\begin{equation}
\DIAG_n\ag{a_1,a_2,\ldots,a_n} = \begin{pmatrix}
a_1 & 0 & \cdots & 0 \\
0 & a_2 & \cdots & 0 \\
\vdots & \vdots & \ddots & \vdots \\
0 & 0 & \cdots & a_n
\end{pmatrix}.
\end{equation}

\paragraph{Set of orthogonal matrices $\Oset{n}$:}
A matrix $\Qbsf\in \Mvspace{n}$ is orthogonal if and only if
\begin{equation}
    \Qbsf^{\Tsf}\Qbsf=\Qbsf\Qbsf^{\Tsf}=\IdentityMatrix_n,
\label{eq:Qproperty}
\end{equation}
where $\IdentityMatrix_n = \DIAG_n\ag{1,\ldots,1}$ is the identity matrix. 
In the previous equation, $\Qbsf^{\Tsf}$ denotes the transpose of $\Qbsf$, i.e. 
if $\Qbsf=\pr{\pr{Q_{ij}}_{j\in \Ic}}_{i\in \Ic}$ then $\Qbsf^{\Tsf}=\pr{\pr{Q_{ji}}_{j\in \Ic}}_{i\in \Ic}$. 
We denote the set of all orthogonal matrices in $\Mvspace{n}$ as $\Oset{n}$.

Let $\qbsf_k=\pr{\pr{Q_{ik}}}_{i\in \Ic} \in \Rbb^n$ be the $k^{\rm th}$ column of $\Qbsf$. Then it follows that
\begin{equation}
\langle \qbsf_i,\qbsf_j\rangle =\delta_{ij},
\label{eq:ColumnsPerpendicular}
\end{equation}
for any $i,j\in \Ic$, 
where $\langle \cdot,\cdot \rangle$ is the standard inner product on $\Rbb^n$ and $\delta_{ij}$ is the Kronecker delta symbol, which equals unity when $i=j$ and zero otherwise. 

\paragraph{Frobenius inner product and norm:}
In this work $\Mvspace{n}$ is a Euclidean vector space.  
Its inner product is $\fprod{\cdot,\cdot} : \Mvspace{n} \times \Mvspace{n}\to \mathbb{R}$,
\begin{equation}
    \fprod{\Xbsf,\Ybsf}= \sum_{i=1}^n \sum_{j=1}^n X_{ij} Y_{ij},
\label{def:FrobInnerProduct}
\end{equation}
where $X_{ij}$, and $Y_{ij}$ are, respectively, the $i$-$j^{\rm th}$ components of $\Xbsf$, and $\Ybsf$. 
The inner product $\fprod{\cdot,\cdot}$ is called the \emph{Frobenius inner product}. 
The Frobenius inner product induces a norm on $\Mvspace{n}$, $\|\cdot\|_{F}: \Mvspace{n}\to \mathbb{R}$, defined as
$$
\| \Absf \|_{F} =\sqrt{\fprod{\Absf,\Absf}},
$$ 
which we call the \emph{Frobenius norm}. 



\subsection{Best approximant of a matrix by the square of a skew-symmetric matrix} 
\label{sec:MainResult}


\begin{proposition}
\label{prop:Existence}
    Problem \ref{p:OptimizationProblem}($\Absf$) admits a solution. 
\end{proposition}

A proof for Proposition \ref{prop:Existence} can be constructed using \textit{Weierstrass' Extreme Value Theorem}~\cite[Proposition A.8]{Dimitri1999} and the facts that: \textit{(i)} there exists an isometric isomorphism between $\Mvspace{n}$ and $\mathbb{R}^{n^2}$,
\textit{(ii)} the set $\Ktwoset{n}$ is closed in $\Mvspace{n}$, and \textit{(iii}) the function which maps $\Ubsf$ to $\| \Ubsf - \Absf \|_F$
 is continuous and coercive on $\Ktwoset{n}$. 

However, Problem \ref{p:OptimizationProblem}($\Absf$) is more effectively addressed by Theorem \ref{thm:main}, which we state later in this section. Not only does Theorem \ref{thm:main} imply Proposition \ref{prop:Existence} as a corollary, but it also provides a recipe for constructing a solution to Problem \ref{p:OptimizationProblem}($\Absf$) . Therefore, we focus the remainder of this section on Theorem \ref{thm:main} and its proof.

To make the statement  of Theorem \ref{thm:main} more compact, we  define what we term the \emph{skew-square-spectral approximant} of a matrix $\Absf$ as follows. 

\begin{definition}[Skew-square-spectral approximant]
\label{def:S3Approx}
Let $\Absf\in \Mvspace{n}$, and $\Bbsf=\pr{\Absf+\Absf^{\Tsf}}/2$. 
It follows from the \emph{Real Spectral Theorem}~\cite[7.29]{Axler2015} that $\Bbsf$ can be decomposed as
\begin{subequations}
    \begin{equation}
\Bbsf = \Nbsf  \Lambdabsf \Nbsf^\Tsf,
\end{equation}
where   $\Nbsf \in \Oset{n}$, and
\begin{equation}
\Lambdabsf:=\DIAG_n\ag{\lambda_1,\lambda_2,\ldots,\lambda_n},
\end{equation} where $\lambda_i \in \mathbb{R}$, $i=1,\ldots,n$, are a non-increasing sequence of real numbers.  
\label{eq:Bsd}
\end{subequations}
The \emph{skew-square-spectral approximant} of $\Absf$ is defined as
\begin{subequations}
\label{eq:Ustevendecom}
\begin{equation}
    \Ubsf^{\star}\ag{\Absf}:=\Nbsf \Dbsf^{\star} \Nbsf^{\Tsf},
\end{equation}
where $\Dbsf^{\star}$ is a diagonal matrix that depends on $\Lambdabsf$. 
Its form is slightly different depending on whether $n$ is even or odd. 
When $n$ is even let $k=n/2$. 
In this case
\begin{equation}
\Dbsf^{\star}:=\DIAG_n\ag{\mu^{\star}_1,\mu^{\star}_1,\ldots,\mu^{\star}_{k},\mu^{\star}_{k}},
\end{equation}
where for $i\in (1,\ldots, k)$
\begin{equation}
\mu^{\star}_i = \left\{ \begin{array}{l l}
\pr{\lambda_{2i-1} + \lambda_{2i}}/2, & \quad \lambda_{2i-1} + \lambda_{2i} \leq 0, \\
0, & \quad \text{otherwise}.
\end{array} \right.
\end{equation}
\end{subequations}
\begin{subequations}
When $n$ is odd let $k=(n-1)/2$. 
In this case
\begin{equation}
    \Dbsf^{\star}:= \DIAG_n\ag{0,\mu^{\star}_1,\mu^{\star}_1,\ldots,\mu^{\star}_k,\mu^{\star}_k},
\end{equation}
where for $i\in (1,\ldots, k)$
\begin{equation}
\mu^{\star}_i = \left\{ \begin{array}{l l}
\pr{\lambda_{2i} + \lambda_{2i+1}}/2, & \quad \lambda_{2i} + \lambda_{2i+1} \leq 0, \\
0, & \quad \text{otherwise}.
\end{array} \right.
\end{equation}
\label{equ:odd}
\end{subequations}
\end{definition}



\begin{restatable}{theorem}{ThmMain}\label{thm:main}
Given $\Absf \in \Mvspace{n}$, a solution to Problem \ref{p:OptimizationProblem}$\pr{\Absf}$ is $\Ubsf^{\star}\ag{\Absf}$. 
\end{restatable}

Recall that problem \ref{p:OptimizationProblem} has been stated in the Introduction. In Theorem~\ref{thm:main}, $\Ubsf^{\star}\ag{\Absf}$ is the skew-square-spectral approximant of $\Absf$, which is defined in Definition~\ref{def:S3Approx}. 
Definition \ref{def:S3Approx} is also a recipe for constructing $\Ubsf^{\star}\ag{\Absf}$. 

In the sequel, we prove Theorem~\ref{thm:main}.

\section{Proofs}
\label{sec:ProofOfMainResults}

\subsection{Supporting results}\label{sub:SupportingResults}

Before proving Theorem~\ref{thm:main} we state four lemmas which we will use in the proof.
Lemma~\ref{thm:frobenius-sym-skew} and Lemma~\ref{thm:frobenius-orthogonal} are standard results, while parts of the proof of Lemma~\ref{thm:K2-eigenvalues} can be found in the literature. For the sake of completeness, we provide complete proofs for these lemmas in \ref{sec:proofs}. The final lemma, Lemma~\ref{thm:lambda-m-n}, is easily proven from existing results in the literature; however, we were unaware of these results during the preparation of this manuscript. A short proof, using specialized results from the literature, as well as our own longer proof, are also provided in \ref{sec:proofs}.

\begin{restatable}{lemma}{ThmFrobeniusSymSkew}
\label{thm:frobenius-sym-skew}
For $\Sbsf \in \Svspace{n}$ and $\Kbsf \in \Kvspace{n}$,
$$
\| \Sbsf + \Kbsf \|_{F}^2 = \| \Sbsf \|_{F}^2 + \| \Kbsf \|_{F}^2.
$$ 
\end{restatable}

\begin{restatable}{lemma}{ThmFrobeniusOrthogonal}
\label{thm:frobenius-orthogonal}
Let $\Absf \in \Mvspace{n}$ and $\Qbsf \in \Oset{n}$. 
Then
$$
  \| \Absf \Qbsf \|_F = \| \Absf \|_F= \| \Qbsf \Absf \|_F.
$$ 
\end{restatable}


For presenting Lemma \ref{thm:K2-eigenvalues}, we first define the set $\Vc_n$. 
\begin{definition}
For $n\in\mathbb{N}$
\begin{equation}
\label{eq:the-set-alt}
\Vc_{n} := \left\{ \Nbsf \Dbsf \Nbsf^\Tsf \SUCHTHAT \Nbsf \in \Oset{n}~\text{and}~\Dbsf \in \Dc_n \right\}.
\end{equation}
The set $\Dc_n$ depends on whether $n$ is even or odd. 
When $n$ is even, $\Dc_n$ is the set of all matrices of the form
$\DIAG_n\ag{\mu_1,\mu_1,\ldots,\mu_k,\mu_k}$, where $\mu_i\leq 0$, $i=1,\ldots, k$, are a non-increasing sequence of real numbers and $k=n/2$. When $n$ is odd $\Dc_n$ is
the set of all matrices of the form
$\DIAG_n\ag{0,\mu_1,\mu_1,\ldots,\mu_k,\mu_k}$, with $\mu_i\leq 0$, $i=1,\ldots, k$, as before and $k=(n-1)/2$.
\end{definition}

\renewcommand{\labelenumii}{\theenumii}
\renewcommand{\theenumii}{$(\theenumi.\alph{enumii})$}

\begin{restatable}{lemma}{ThmKTwoEigenvalues}
\label{thm:K2-eigenvalues}
The sets $\Vc_n$ and $\Ktwoset{n}$ are the same. 
\end{restatable}

\begin{restatable}{lemma}{ThmLambdaMN}\label{thm:lambda-m-n}
Let $\Lambdabsf = \DIAG_n\ag{\lambda_1,\lambda_2,\ldots,\lambda_n}$ where $\lambda_1 \geq \lambda_2 \geq \ldots \geq \lambda_n$, and let $\Dbsf = \DIAG_n\ag{\mu_1,\mu_2,\ldots,\mu_n}$ where $\mu_1 \geq \mu_2 \geq \ldots \geq \mu_n$. Then,
\begin{equation}\label{eq:minp}
\min_{\Qbsf \in \Oset{n}} \| \Lambdabsf - \Qbsf \Dbsf \Qbsf^\Tsf \|_F = \| \Lambdabsf - \Dbsf \|_F.
\end{equation}
\end{restatable}

\subsection{Proof of Theorem~\ref{thm:main}}\label{sub:ProofOfMainTheorem}

\ThmMain*
\begin{proof}
We consider the case where $n$ is even, i.e., $n=2k$, for some $k\in\mathbb{N}$. 
The case where $n$ is odd can be handled similarly. 
Let $\Cbsf = (\Absf - \Absf^\Tsf)/2$ and $\Bbsf = (\Absf + \Absf^\Tsf)/2$. For any $\usf{U} \in \Ktwoset{n}$, by Lemma~\ref{thm:frobenius-sym-skew}
\begin{equation}
\| \Absf - \usf{U} \|_{F}^2 = \| \Bbsf - \usf{U} \|_{F}^2 + \| \Cbsf \|_{F}^2.
\label{eq:fp1}
\end{equation}
Owing to \eqref{eq:fp1}, a solution to \ref{p:OptimizationProblem}$(\Absf)$ is also a solution to \ref{p:OptimizationProblem}$(\Bbsf)$ and vice versa.
Furthermore, from $\Ubsf^{\star}\ag{\Absf}$'s definition we note that $\Ubsf^{\star}\ag{\Absf}=\Ubsf^{\star}\ag{\Bbsf}$. Thus, it suffices to show that
a solution to \ref{p:OptimizationProblem}$\pr{\Bbsf}$ is $\Ubsf^{\star}\ag{\Bbsf}$, i.e.,
\begin{equation}
\norm{\Bbsf - \usf{U}^{\star}\ag{\Bbsf}}_{F}^2\leq\norm{\Bbsf - \usf{U}}_{F}^2,
\label{eq:result}
\end{equation}
for any $\usf{U} \in \Ktwoset{n}$.

\begin{enumerate}
    \item Since $\Ubsf\in \Ktwoset{n}$ it follows from
Lemma~\ref{thm:K2-eigenvalues} that $\usf{U}$ can be decomposed as
\begin{equation}
\usf{U}= \Mbsf \Dbsf \Mbsf^\Tsf,
\label{eq:Udecom}
\end{equation}
where $\Dbsf = \DIAG_n\ag{\mu_1,\mu_1,\ldots,\mu_k,\mu_k}$ with $0 \geq \mu_1 \geq \ldots \geq \mu_k$, and $\Mbsf \in \Oset{n}$.

\begin{subequations}
\item Applying $\Ubsf$'s decomposition given in \eqref{eq:Udecom} and $\Bbsf$'s decomposition given in \eqref{eq:Bsd} we get that
\begin{align}
\norm{ \Bbsf - \usf{U} }_{F}^2&=\norm{ \Nbsf \Lambdabsf \Nbsf^\Tsf - \Mbsf \Dbsf \Mbsf^\Tsf }_{F}^2,\notag\\
&=\norm{\Nbsf\pr{\Lambdabsf - \Nbsf^\Tsf \Mbsf \Dbsf \Mbsf^\Tsf \Nbsf} \Nbsf^\Tsf }_F^2, \notag \\
&= \norm{\Lambdabsf - \Nbsf^\Tsf \Mbsf \Dbsf \Mbsf^\Tsf \Nbsf }_F^2,
\end{align}
where the third equality follows from two applications of Lemma~\ref{thm:frobenius-orthogonal} to remove the orthogonal matrices $\Nbsf$ and $\Nbsf^\Tsf$. Defining $\Qbsf = \Nbsf^\Tsf \Mbsf \in \Oset{n}$,
\begin{align}
\norm{ \Bbsf - \usf{U} }_{F}^2
&=\norm{\Lambdabsf  - \Qbsf \Dbsf \Qbsf^\Tsf}_F^2,\label{eq:Term13}\\
&\ge \norm{\Lambdabsf  -  \Dbsf }_F^2,
\label{eq:Term14}
\end{align}
where the inequality follows from Lemma~\ref{thm:lambda-m-n}. 
\end{subequations}

\item Next, using $\Bbsf$'s decomposition given in \eqref{eq:Bsd} and $\Ubsf^{\star}\ag{\Bbsf}$'s decomposition (which is the same as $\Ubsf^{\star}\ag{\Absf}$'s decomposition) given in \eqref{eq:Ustevendecom} we have
\begin{equation}
\norm{\Bbsf - \usf{U}^{\star}\ag{\Bbsf} }_{F}^2=
\norm{ \Lambdabsf - \Dbsf^{\star} }_{F}^2,
\label{eq:Term21}
\end{equation}
where we have used similar manipulations to the previous step.

\item Subtracting \eqref{eq:Term21} from \eqref{eq:Term14} we get
\begin{align}
\| \Bbsf - \usf{U} \|_{F}^2 -\| \Bbsf - \usf{U}^{\star}\ag{\Bbsf} \|_{F}^2&\geq \|  \Lambdabsf -  \Dbsf \|_{F}^2- \| \Lambdabsf - \Dbsf^{\star}\|_{F}^2,\notag\\
&= 2\sum_{i=1}^k \pr{\mu_i-\mu_i^{\star}}(\mu_i+\mu_i^{\star} \notag\\
&~~~~~~~-\lambda_{2i-1} - \lambda_{2i}),
\label{eq:T1T2}
\end{align}
where the equality follows direct calculation of the Frobenius norms using the explicit formulas for $\Lambdabsf$, $\Dbsf$, and $\Dbsf^{\star}$ in terms of $\lambda_i$, $\mu_i$, and $\mu_i^{\star}$, respectively. 

\item We claim each term in the sum in \eqref{eq:T1T2} is non-negative, and so \eqref{eq:result} holds. Let $1 \leq i \leq k$.
\begin{enumerate}
    \item Suppose that $\lambda_{2i-1}+\lambda_{2i}$ is non-positive. Then, by definition $\mu_i^{\star}=\pr{\lambda_{2i-1}+\lambda_{2i}}/2$, and so the term in the sum simplifies to
    $$
    \pr{\mu_i-\frac{\lambda_{2i-1}+\lambda_{2i}}{2}}^2,
    $$
    which is clearly non-negative.
    \item Otherwise, suppose that $\lambda_{2i-1}+\lambda_{2i}$ is positive. Then $\mu_i^\star = 0$ and the term in the sum is
    $$
    \mu_i \pr{\mu_i - \lambda_{2i-1} - \lambda_{2i}}.
    $$
    Since $\mu_i \leq 0$ (see Step~1) and $- \lambda_{2i-1} - \lambda_{2i} < 0$, this is the product of two negative numbers, which is positive.
\end{enumerate}
\end{enumerate}
\end{proof}

\section{Algorithm and examples}
\label{sec:exaapp}

\subsection{Algorithm}
\label{subsec:algorithm}
Following Definition~\ref{def:S3Approx}, we summarize the procedure to compute the best approximant of $\Absf$ in $\Ktwoset{n}$, i.e., $\Ubsf^{\star}\ag{\Absf}$ in Algorithm~\ref{algo:ustar}.

\RestyleAlgo{ruled}
\SetKwComment{Comment}{/* }{ */}
\begin{algorithm}[h]
\caption{Computing a best approximant to a matrix $\Absf \in \Mvspace{n}$ in the set $\Ktwoset{n}$}
\label{algo:ustar}
\textbf{Input}: A matrix $\Absf \in \Mvspace{n}$, $n\geq1$\;
\textbf{Output}: A best approximant to $\Absf$ in $\Ktwoset{n}$, $\Ubsf^{\star}\ag{\Absf}$ \;

Compute the symmetric part of $\Absf$, $\Bbsf=\pr{\Absf+\Absf^{\Tsf}}/2$\;

Compute a real spectral decomposition of
$$
\Bbsf = \Nbsf  \DIAG_n\ag{\lambda_1,\lambda_2,\ldots,\lambda_n} \Nbsf^\Tsf
$$
in which $\Nbsf \in \Oset{n}$ and $\lambda_i \in \mathbb{R}$, $i=1,\ldots,n$, are a non-increasing sequence of real numbers\;

$\usf{d}^{\star}\leftarrow \emptyset$ \Comment*[r]{an empty list}
\eIf{$n$ is even}{
$k=n/2$\;
\For{$i=1$ to $k$}{
\eIf{$\lambda_{2i-1}+\lambda_{2i}\leq0$}{$\mu^{\star}\leftarrow \pr{\lambda_{2i-1}+\lambda_{2i}}/2$\;}{$\mu^{\star}\leftarrow 0$\;}
Append $\mu^{\star}$ to $\usf{d}^{\star}$ twice;
}
}
{
$k=(n-1)/2$\;
Append $0$ to $\usf{d}^{\star}$\;
\For{$i=1$ to $k$}{
\eIf{$\lambda_{2i}+\lambda_{2i+1}\leq0$}{$\mu^{\star}\leftarrow \pr{\lambda_{2i}+\lambda_{2i+1}}/2$\;}{$\mu^{\star}\leftarrow 0$\;}
Append $\mu^{\star}$ to $\usf{d}^{\star}$ twice\;
}
}
$\Dbsf^{\star}\leftarrow \DIAG_n\ag{\usf{d}^{\star}}$\;
\Return{$\Ubsf^{\star}\ag{\Absf}=\Nbsf\Dbsf^{\star}\Nbsf^\Tsf$.}
\end{algorithm}

\subsection{Basic examples}
\label{subsec:examples}

We illustrate the application of Algorithm~\ref{algo:ustar} for some simple matrices, which highlight the salient steps and emphasize the potential existence of multiple solutions.

\paragraph{Example 1}
Let
    \begin{equation}
      \Absf = \begin{pmatrix} -1 & 4 & 2 \\ 2 & -1 & 3 \\ -2 & -3 & -6 \end{pmatrix}.
      \label{eq:Ex1A}
    \end{equation}
The symmetric part of $\Absf$ is
    \begin{equation}
    \Bbsf= \pr{\Absf+\Absf^{\Tsf}}/2=
    \begin{pmatrix} -1 & 3 & 0 \\ 3 & -1 & 0 \\ 0 & 0 & -6 \end{pmatrix}.
    \label{eq:Ex1B}
    \end{equation}
A real spectral decomposition of $\Bbsf$ of the form stipulated by Algorithm~\ref{algo:ustar}\footnote{\label{ft:orsd}It can be shown that this real spectral decomposition  for the symmetric part of the given $\Absf$ is the only one possible for it that has the form stipulated by Algorithm~\ref{algo:ustar}.} is
\begin{subequations}
\begin{equation}
    \Bbsf =    \Nbsf \Lambdabsf
    \Nbsf^\Tsf,
\end{equation}
where
\begin{align}
\Nbsf
&= \begin{pmatrix} \frac{1}{\sqrt{2}} & -\frac{1}{\sqrt{2}} & 0 \\ \frac{1}{\sqrt{2}} & \frac{1}{\sqrt{2}} & 0 \\ 0 & 0 & 1 \end{pmatrix},\\
 \Lambdabsf&=   \DIAG_3\ag{2,-4,-6}.
\end{align}
\label{eq:DecomBEx1}
\end{subequations}
Next we need to construct the matrix $\Dbsf^{\star}$.
As $n=3$ (and $k=(n-1)/2=1$) we know
\begin{equation}
\Dbsf^\star=\DIAG_3\ag{0,\mu_1^{\star},\mu_1^{\star}}.
\label{eq:Dform_for_n_1}
\end{equation}
Since $\lambda_2+\lambda_3=-4-6=-10<0$, we get that
$$
\mu_1^{\star}=(-4-6)/2=-5.
$$
Thus, we can construct $ \Ubsf^\star\ag{\Absf}$:
\begin{equation}
\Ubsf^\star\ag{\Absf} = \Nbsf \Dbsf^\star \Nbsf^\Tsf = \begin{pmatrix}
        -\frac{5}{2} & \frac{5}{2} & 0 \\ \frac{5}{2} & -\frac{5}{2} & 0 \\ 0 & 0 & -5
    \end{pmatrix}.
\end{equation}


\paragraph{Example 2} Our work provides a means of constructing \emph{a} best approximation to a given matrix in the sense of Problem \ref{p:OptimizationProblem}($\Absf$); however, there is no guarantee of uniqueness. For example, consider the matrix $\Absf=\Bbsf = \DIAG_3\ag{-1,-1,-1}$. This matrix does not have a unique spectral decomposition of the form stipulated by Algorithm~\ref{algo:ustar}. That is, $\Bbsf$ can be written as $\Nbsf \DIAG_3\ag{-1,-1,-1}\Nbsf^{\Tsf}$ where $\Nbsf$ is any element in $\Oset{3}$. Following Algorithm~\ref{algo:ustar}, we set
\begin{equation}
\Ubsf^\star\ag{\Absf} = \Nbsf \DIAG_3\ag{0,-1,-1} \Nbsf^\Tsf,
\end{equation}
where we reiterate that $\Nbsf$ is an arbitrary element of $\Oset{3}$. In other words, all elements of the set
\begin{equation}
\left\{    \Nbsf\, \DIAG_3\ag{0,-1,-1}\, \Nbsf^{\Tsf}\SUCHTHAT \Nbsf \in \Oset{3}\right\}\subset \Ktwoset{n}
\end{equation}
are best approximants of $\Absf$.

\section{An application in measuring angular velocity}
\label{sec:expdemo}

\begin{figure}[h]
\centering
\includegraphics[width=17cm]{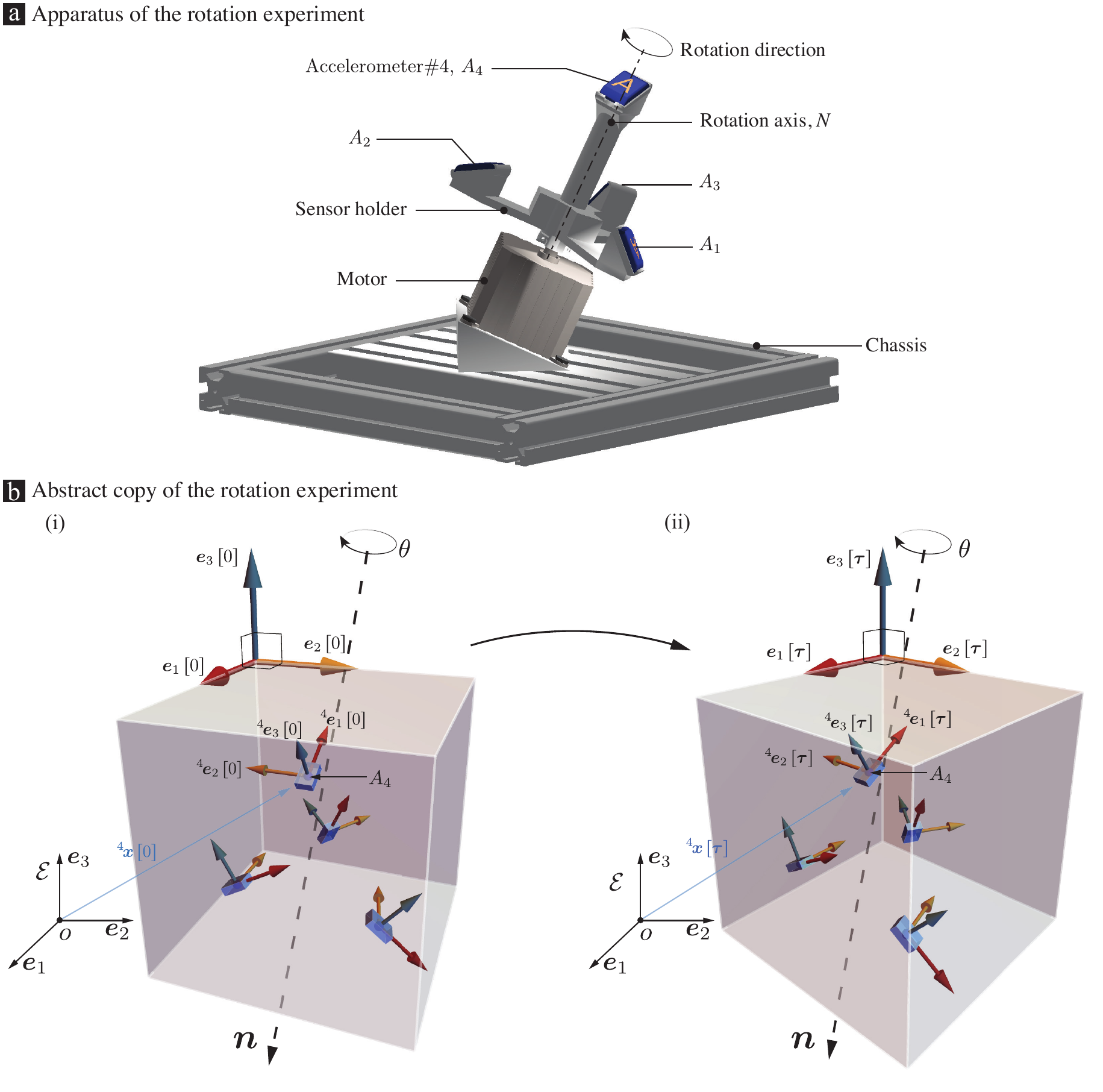}
    \caption{Rigid body rotation experiment setup. (a) shows a sketch of the rigid body rotation experiment described in \S\ref{sec:rotationexp}. 
    (b) shows a more abstract version of the experiment. 
    Subfigures (i) and (ii) respectively denote the configurations of the rigid body in our experiment at the initial and a later time instance. 
    All the mathematical symbols in this figure are described in \S\ref{sec:rotationexp}--\ref{sec:ApplyingAO}. 
    For instance, the arrows marked $\u{e}_i$, $i=1,2,3$, denote the basis vectors, and
    the blue cuboid marked $A_4$ denotes the accelerometer \#4. 
    Its position vector at the initial  and at a later time instance is, respectively, marked in (i) and (ii) as  ${}^4\!\u{x}\ag{0}$ and ${}^4\!\u{x}\ag{\bm{\tau}}$. 
    The components of ${}^4\!\u{x}\ag{0}$ w.r.t.~$\pr{\u{e}_i}_{i\in I}$, which we denote as ${}^{4}\!\usf{x}\ag{0}$, in the punctuated rotation trial (see, e.g., \S\ref{sec:experiment}) are $\pr{0.02,0.03,0.11}$. 
    Similarly, the components of the other accelerometers' initial position vectors are  ${}^{1}\!\usf{x}\ag{0}=\pr{-0.08,-0.01,0.04}$, ${}^2\!\usf{x}\ag{0}=\pr{0.04,-0.06,0.01}$, and ${}^3\!\usf{x}\ag{0}=\pr{0.02,0.08,-0.05}$. 
    The arrows marked $\lsc{4}\u{e}_i\ag{0}$, $i=1,2,3$, in (i) denote accelerometer \#4's measurement directions at the initial time instance. 
    Their components w.r.t. $\pr{\u{e}_i}_{i\in I}$ in the punctuated rotation trial are $\pr{-0.68,0.14,0.71}$, $\pr{0.08,-0.95,0.27}$, and $\pr{0.72,0.24,0.64}$. 
    Similarly, the components of $\pr{\lsc{1}\u{e}_i\ag{0}}_{i\in I}$ are  $\pr{\pr{0.44,0.66,0.61},\pr{-0.87,0.48,0.11},\pr{-0.22,-0.58,0.79}}$, of $\pr{\lsc{2}\u{e}_i\ag{0}}_{i\in I}$ are $\pr{\pr{0.07,0.97,0.22},\pr{-0.96,0,0.28},\pr{0.27,-0.23,0.93}}$, and of $\pr{\lsc{3}\u{e}_i\ag{0}}_{i\in I}$ are $\pr{\pr{-0.42,0.45,-0.78},\pr{-0.90,-0.24,0.34},\pr{-0.03,0.85,0.51}}$. 
     Furthermore, the arrows marked $\u{e}_i\ag{0}$, $i=1, 2, 3$, in (i) denote the body basis vectors at the initial time instance. 
     Their respective components in the punctuated rotation trial are $\pr{1,0,0}$, $\pr{0,1,0}$, and $\pr{0,0,1}$. 
     The arrow marked $\u{n}$ denotes the rotation axis, and its components are $\usf{n}=\pr{-0.27,-0.28,-0.92}$. 
     } 
 \label{fig:testsetup}
\end{figure}

In this section, we  demonstrate an experimental application of Theorem~\ref{thm:main}.
We use our result to robustly estimate a rigid body's angular velocity from experimental accelerometer data.
Specifically, we apply the $\sqrt{\text{AO}}$-algorithm to the accelerometer data from a rigid body rotation experiment (Fig.~\ref{fig:testsetup}) and estimate the angular velocity vector of the rigid body in that experiment.
Recall that Theorem~\ref{thm:main} is an integral part of the $\sqrt{\rm AO}$-algorithm; it is used for carrying out the critical step of constructing an approximation for the matrix $\widetilde{\Bbsf}\ag{\tau_i}$ in the space $\Ktwoset{n}$ (see \S\ref{sec:MotivationApplication} for details).
Without this approximation step, the $\sqrt{\rm AO}$-algorithm reduces to the plain $\sqrt{\rm AO}$-algorithm. The $\sqrt{\rm AO}$-algorithm is a viable means for estimating  angular velocity in real-world situations wherein the data contains noise and errors, and, due to the elasticity of  materials, the accelerometers are no longer rigidly affixed to one another (so that generally $\widetilde{\Bbsf}\ag{\tau_i} \not \in \Ktwoset{n}$).
In contrast, the plain $\sqrt{\rm AO}$-algorithm can only calculate the angular velocity in ideal situations\footnote{By ``real-world'' and ``ideal'' we are not simply contrasting between experimental and computational situations. Even motion extracted from computational mechanics simulations of rigid body motion can contain numerical noise, especially at large time steps, and it too could only be an approximation of rigid body motion, since  specialized time integration schemes are required to perfectly maintain the rigidity constraint in the simulation.}, wherein the data is free from noise and errors and the motion is that of a perfect rigid body.

\clearpage

\subsection{Rigid body motion and rotation experiments}
\label{sec:rotationexp}
We conducted a rigid body rotation experiment to demonstrate an application of the $\sqrt{\rm AO}$-algorithm, and, consequently, that of  Theorem~\ref{thm:main}. 

\subsubsection{Rigid body motion}
\label{sec:rbm}
A sketch of the rigid body and the apparatus we used to impose a motion on it in our rotation experiment are shown in Fig.~\ref{fig:testsetup}(a). 
Photographs of our experimental setup are provided in Fig.~\ref{fig:rotationtest}. 
A more abstract representation of our rigid body is shown in Fig.~\ref{fig:testsetup}(b). 
As shown in Fig.~\ref{fig:testsetup}(b), we take our rigid body  to execute its motion in the physical Euclidean point space $\mathcal{E}$, which is a three dimensional real affine space. We term the inner product space associated with it the physical Euclidean vector space, and denote it as $\mathbb{E}$. 
An orthonormal basis for $\mathbb{E}$ is the set of vectors $\pr{\boldsymbol{e}_i}_{i\in I}$, where $I=(1,2,3)$, shown in Fig.~\ref{fig:testsetup}(b). 
In our formulation the vectors $\boldsymbol{e}_i$ have units of $\textsf{meters}$\footnote{This formalism where the elements in a vector space have units, while their components w.r.t.~any basis of that vector space are non-dimensional, was introduced in \cite{wan2024mechanics,Wan2023brain}.
}. 
Let the position vector of a rigid body material particle $\mathcal{P}$ w.r.t. the origin $o$, shown in Fig.~\ref{fig:testsetup}(b), at the time instance $\boldsymbol{\tau}$ be $^{\mathcal{P}}\!\boldsymbol{x}\ag{\boldsymbol{\tau}}$. 
We then refer to
\begin{equation}
^{\mathcal{P}}\!\xbsf\ag{\tau}:=\pr{
^{\mathcal{P}}\!\boldsymbol{x}\ag{\boldsymbol{\tau}}\cdot \boldsymbol{e}_i}_{i\in I}\in \mathbb{R}^3,
\label{eq:DimenNonDimen}
\end{equation}
as particle $\mathcal{P}'s$ non-dimensional position vector
at the non-dimensional time instance $\tau\in \mathbb{R}$. 
In our formulation $\tau$ is defined such that  $\tau\ \textsf{seconds}=\boldsymbol{\tau}$. 
From here on we will drop the qualifiers ``dimensional'' and ``non-dimensional,'' since these should  be clear from context\footnote{To aid the reader, we denote non-dimensional quantity using sans serif fonts. For example, we denote the non-dimensional counterpart to the dimensional/physical scalar $\boldsymbol{\tau}$ as $\tau$, and the non-dimensional counterpart to the dimensional vector $\boldsymbol{x}$ as $\xbsf$.}. 
In a general rigid body motion $^{\mathcal{P}}\!\xbsf\ag{\tau}$ evolves as
\begin{align}
\lsc{ \sc \mathcal{P} } \! \xbsf\ag{\tau} &= \Rbsf\ag{ \tau } \, \lsc{\mathcal{P}} \! \xbsf_0 + \cbsf\ag{\tau},
\label{eq:GRBM}
\end{align}
where $\Rbsf\ag{\tau}$ is a rotation matrix, belonging to the special orthogonal group $SO(3,\mathbb{R})$ (a subset of $\Oset{3}$), and $\lsc{\mathcal{P}} \! \xbsf_0, \cbsf\ag{\tau}\in \mathbb{R}^3$. The matrix $\Rbsf\ag{\tau}$ and the vector $\cbsf\ag{\tau}$ are called  the rotation matrix and the translation vector at the time instance $\tau$, respectively. 

The  angular velocity vector of a rigid body from classical physics, $\boldsymbol{\omega}[\u{\tau}]$, is related to $\Rbsf[\tau]$ as follows. 
The \textit{angular velocity matrix in the body frame} $\overline{\Wbsf}\ag{\tau}$ is defined as
\begin{equation}
\overline{\Wbsf}\ag{\tau}=\Rbsf^{\rm \sf T}\ag{\tau}\Rbsf'\ag{\tau},
\label{eq:Wbf}
\end{equation}
where $\Rbsf'\ag{\tau}$ is the derivative of  $\Rbsf[\tau]$. 
It can be shown that $\overline{\Wbsf}\ag{\tau}\in \Kvspace{3} \cong \mathfrak{so}(3,\mathbb{R})$, the Lie algebra of $SO(3,\mathbb{R})$. 
The \textit{body angular velocity vector} is defined as
\begin{equation}
\busf{w}\ag{\tau}=\star\overline{\Wbsf}\ag{\tau},
\label{eq:wbf}
\end{equation}
where  ``$\star$'' denotes the map $\mathfrak{so}\pr{3,\mathbb{R}}\ni \Wbsf\mapsto \star\Wbsf\in \mathbb{R}^3$,
\begin{align*}
    \star\Wbsf=\pr{\Wbsf_{\cdot 3 \cdot 2},\Wbsf_{\cdot 1 \cdot 3},\Wbsf_{\cdot 2 \cdot 1}}.
\end{align*} 
The \textit{angular velocity vector in the laboratory frame} is
\begin{align}
	\wbsf[\tau]&:=\sum_{i\in I}\overline{\wbsf}_{\cdot i}\ag{\tau}\ebsf_i[\tau], \label{eq:wlf}
	\intertext{where}
	\ebsf_i\ag{\tau}&:=\Rbsf\ag{\tau}\ebsf_i, \label{eq:bfevolution}
\end{align}
and    $\ebsf_1:=(1,0,0)$, $\ebsf_2:=(0,1,0)$, and $\ebsf_3:=(0,0,1)$\footnote{The vectors $\ebsf_i\ag{\tau}:=\Rbsf\ag{\tau}\ebsf_i$ are called the body basis vectors. Their dimensional counterparts appear  in Fig.~\ref{fig:testsetup}(b) as   $\boldsymbol{e}_i\ag{\boldsymbol{\tau}}$.
These vectors can be thought of as being attached to the body and rotating with it (cf. Figs.~\ref{fig:testsetup}(b.i) and \ref{fig:testsetup}(b.ii)).}. 
The vector $\wbsf[\tau]$ is, in fact, a non-dimensional form of $\boldsymbol{\omega}\ag{\bm{\tau}}$.

\subsubsection{Experimental configurations}
\label{sec:experiment}
The $\sqrt{\rm AO}$-algorithm provides an approach to estimate the body angular velocity matrix, $\overline{\Wbsf}\ag{\tau}$, of a rigid body when provided with acceleration measurements at four of its points (see \S\ref{sec:MotivationApplication} for details). 
It can handle data from the general rigid body motion given in \eqref{eq:GRBM}. 
However, in order to make it easy for us to quantify the accuracy of the $\sqrt{\rm AO}$-algorithm's estimated angular velocities, in our experiments we  only impose a relatively simple, special case of the motion given in \eqref{eq:GRBM}. 
Specifically, in our experiments we impose
\begin{subequations}
\label{eq:Cat1RBM}
    \begin{align}
    \cbsf\ag{\tau}&\equiv \u{0},
    \intertext{i.e., no translation, and}
    \Rbsf[\tau]&=\nbsf\otimes\nbsf +\pr{\Ibsf-\nbsf\otimes \nbsf}\cos \ag{\theta \ag{\tau}}+\pr{\ast \nbsf}\sin\ag{\theta \ag{\tau}},
\label{eq:SRBM}
\end{align}
where $\Ibsf$ is the identity element in $\Mvspace{3}$, $\nbsf\in \mathbb{R}^3$,  $\theta\ag{\cdot}$ is a continuous function, and ``$\otimes$'' is the dyadic product, such that if $\ubsf,\ \vbsf\in \mathbb{R}^n$ then $\ubsf\otimes \vbsf \in \mathcal{M}_n$, and $\pr{\ubsf\otimes \vbsf}_{\cdot i\cdot j}=\ubsf_{\cdot i}\vbsf_{\cdot j}$. 
The ``$\ast$'' appearing in \eqref{eq:SRBM} is the inverse of $\star$, mapping $\nbsf \in \Rbb^3$ to
\begin{align*}
    \ast\nbsf=\pr{\pr{\sum_{k\in I} \epsilon_{ikj}\nbsf_{\cdot k}}_{j\in I}}_{i \in I} \in \mathfrak{so}\pr{3,\Rbb},
\end{align*}
where $\epsilon_{ijk}$ is the Levi-Cevita symbol. 
The vector $\nbsf$ is called the rotation axis, and $\theta\ag{\tau}\in \mathbb{R}$  the rotation angle at the time instance $\tau$.
%

As per \eqref{eq:SRBM} our rigid body rotates about the fixed axis $\nbsf$. 
(The dimensional counterpart of  $\nbsf$, namely $\bm{n}$, can be seen in, e.g., Fig.~\ref{fig:testsetup}(b).) 
We use a brushless direct current (BLDC) motor (CPM-SCHP-3426D-ELNB, ClearPath, Teknic, USA; see Fig.~\ref{fig:testsetup}(a)) to vary the rotation angle over time in a precise and systematic fashion. 
Configurations of our rigid body at several $\theta$ values are shown in Fig.~\ref{fig:frames}(b). 

We carried out three experiments, each with a different imposed $\theta\ag{\cdot}$.
The results from all three experiments are quite similar. 
Therefore, here we only present results from the first experiment, which we term \textit{punctuated rotation}; and discuss the other two (which we refer to as \textit{oscillatory rotation}, and \textit{constant rotation}) in \S\ref{sec:addrotationtest}. 
In the punctuated rotation experiment
\begin{equation}
    \theta[\tau]=  \frac{\omega_{\rm m}}{2} \tau -\frac{2\omega_{\rm m} \tau_1}{\pi^3} \sum_{n=1,3,5}
    \frac{1}{n^3} \sin\ag{\frac{2n \pi \tau}{\tau_1}}, \label{eq:PunctuatedRotationTheta}
\end{equation}
 \end{subequations}
where $\omega_{\rm m}$ and $\tau_1$ are positive real numbers. 
The angular velocity in the punctuated rotation motion \eqref{eq:Cat1RBM} has the form of a triangle wave with period $\tau_1 \ \textsf{sec}$, 
varying between zero and $\omega_{\rm m}\ \textsf{Hertz}$\footnote{It is more common to use $\textsf{rad/sec}$ as the units of angular velocity. We, however, chose to use the equivalent units of $\textsf{Hertz}$.}, see Fig.~\ref{fig:frames}(a.ii). 

\subsection{Acceleration measurements}
\label{sec:accemeasurement}

\begin{figure*}[t!]
\centering
\includegraphics[width=17cm]
{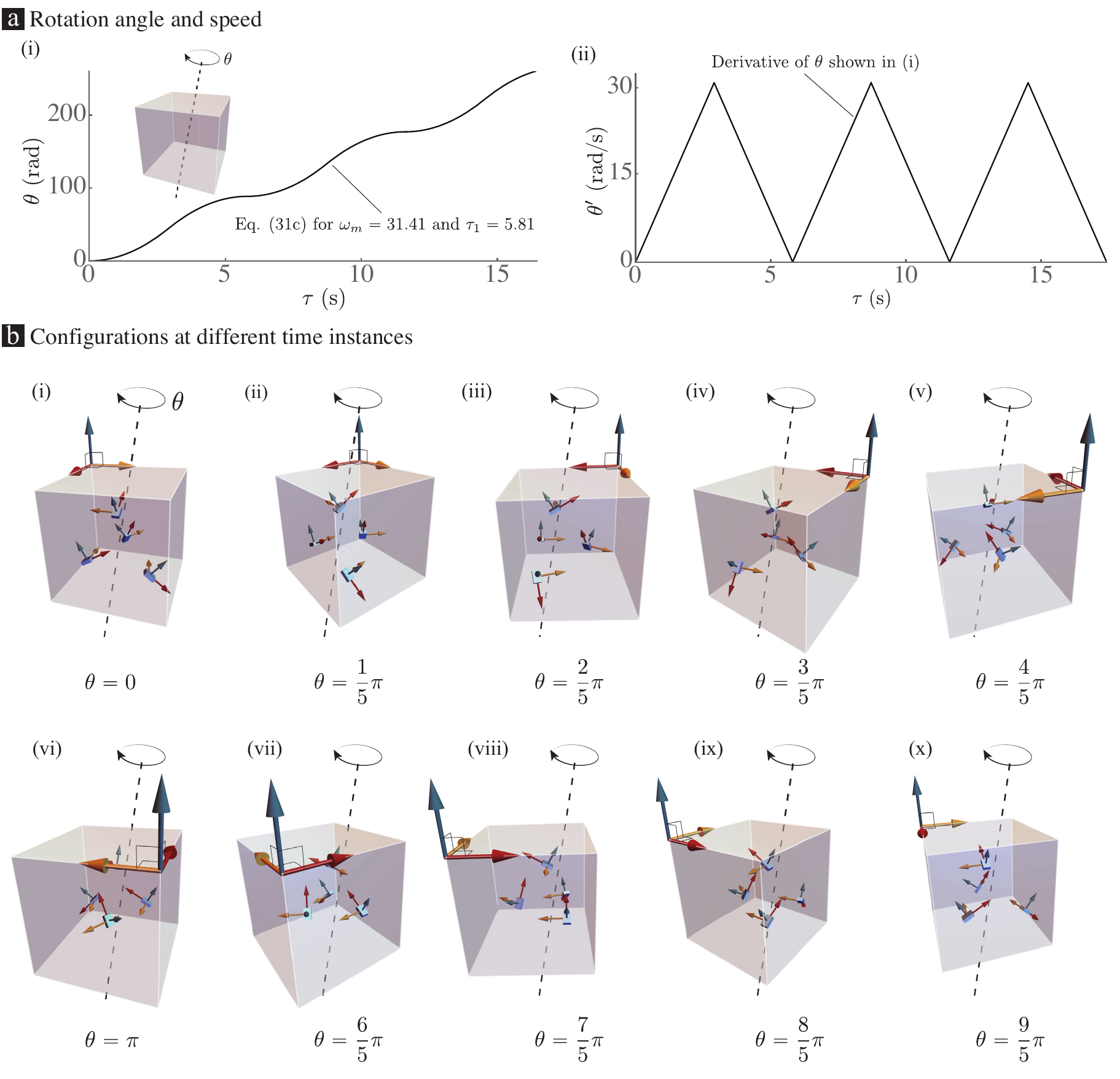}
    \caption{Punctuated rotation experiment. 
	Subfigures (a.i) and (a.ii), respectively, show the variation of the rotation angle and its derivative over time in the punctuated rotation experiment (see \S\ref{sec:experiment} for details). 
	 Subfigure (b) shows the configurations of the rigid body (gray cuboid) and the accelerometers (blue cuboids within the gray cuboid) attached to it for different rotation angles. 
	 The three arrows attached to the top corner of the gray cuboid denote the body basis vectors. 
	 The three arrows attached to a blue cuboid denote the accelerometer's measurement directions. 
	 } 
    \label{fig:frames}
\end{figure*}

\begin{figure*}[t!]
\centering
\includegraphics[width=17cm]{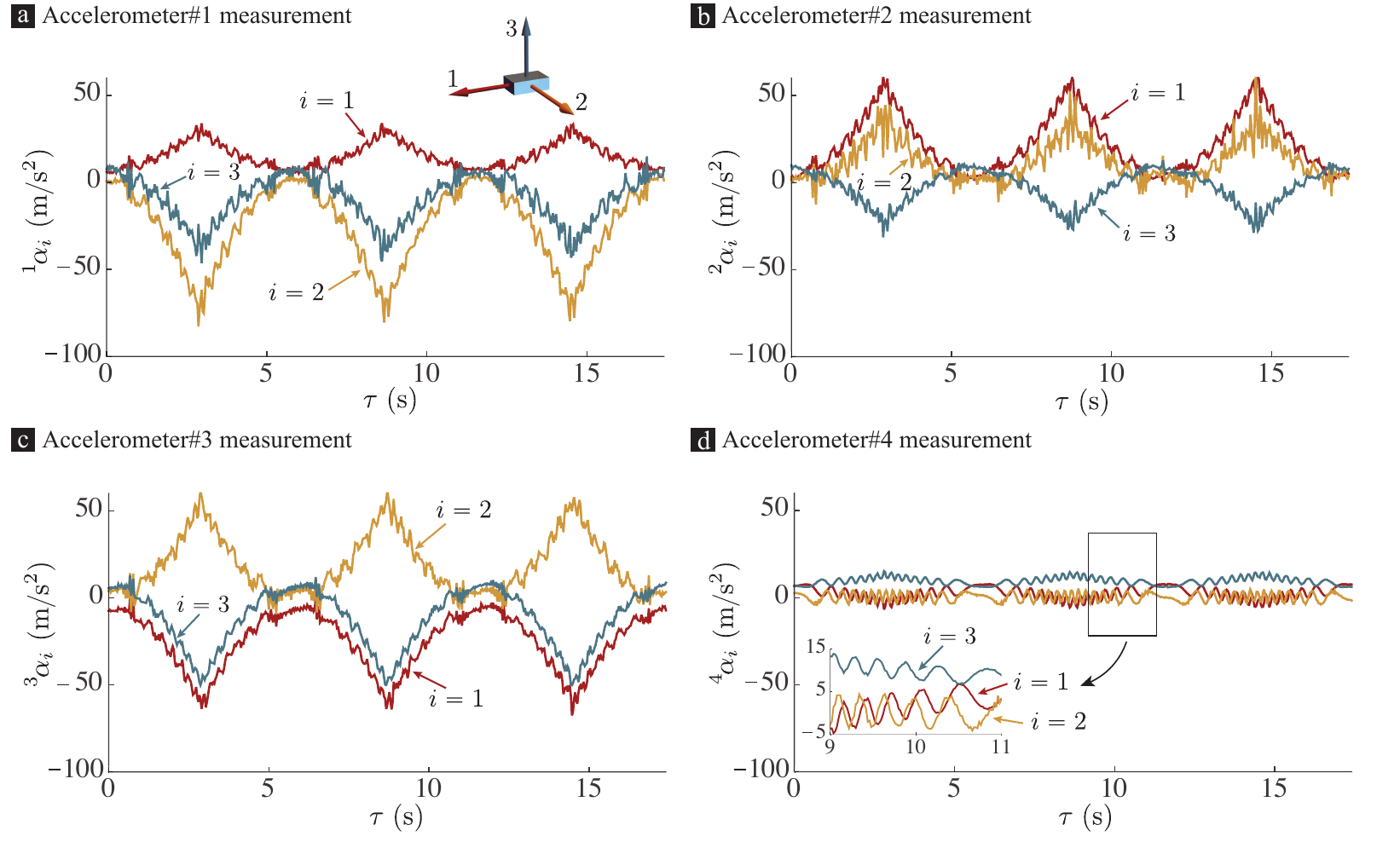}
    \caption{Acceleration measurements in the punctuated rotation experiment (see \S\ref{sec:accemeasurement} for details). 
	Subfigures (a)--(d) correspond to the measurements from accelerometers $A_1$--$A_4$, respectively. 
	The positions and orientations of the four accelerometers w.r.t. the rigid body is illustrated in Figs.~\ref{fig:testsetup}(b) and \ref{fig:frames}(b). 
	Each subfigure shows three graphs, labeled $i=1$, $i=2$, and $i=3$. 
	These, respectively, correspond to the acceleration components measured by the accelerometer in its  $1$, $2$, and $3$, directions 
	(see, e.g., the three arrows attached to the cuboid in the inset in (a)).} 
 \label{fig:measurment}
\end{figure*}

Results from a trial of the punctuated rotation experiment are shown in Figs.~\ref{fig:measurment}--\ref{fig:motorcom}. 
The rotation $\theta\ag{\cdot}$ was that in \eqref{eq:PunctuatedRotationTheta} with $\omega_{\rm m}= 31.41$, and $\tau_1=5.81$. 
(Graphs of $\theta\ag{\cdot}$ and its derivative for these parameter values are shown in Figs.~\ref{fig:frames}(a.i) and (a.ii), respectively.) 
The axis vector $\nbsf$ in this trial was $\pr{-0.27,-0.28,-0.92}$. 
(The dimensional form $\bm{n}$ is shown in Fig.~\ref{fig:testsetup}(b).) 

It follows from \eqref{eq:Wbf}, \eqref{eq:wbf}, and \eqref{eq:SRBM}  that the body angular velocity vector $\overline{\wbsf}$ is
\begin{equation}
	\overline{\wbsf}\ag{\tau}=\theta'\ag{\tau}\, \nbsf. 
\label{eq:wbfIdeal}
\end{equation}
The time-dependency of the components of $\overline{\wbsf}$'s for the punctuated rotation experiment are shown in Fig.~\ref{fig:motorcom}.   
The $\sqrt{\rm AO}$-algorithm can estimate these components when provided with acceleration measurements at four of the rigid body's material particles. 
We discuss these measurements in the reminder of this section, and the application of the $\sqrt{\rm AO}$-algorithm to these measurements in the next section. 
(For prediction of mild traumatic brain injury, the angular velocity vector $\wbsf$'s magnitude is more critical than its direction. 
The body angular velocity vector $\overline{\wbsf}$'s magnitude is always equal to that of $\wbsf$'s. 
In the current trial, due to the simple nature of the imposed motion, their directions are same as well.) 

Acceleration measurements at four of the rigid body's particles in the trial are shown in Fig.~\ref{fig:measurment}. 
We denote the four particles as $ A_1,\ldots,A_4$ (marked, e.g., in Fig.~\ref{fig:testsetup}(a)). 
The acceleration vector at each $A_{\ell}$, $\ell=1,2,3,4$, was measured by a tri-axial accelerometer (Blue Trident, Vicon, UK; 1.6 kHz sampling rate) that was rigidly affixed to it.  
(These accelerometers are the blue objects in Fig.~\ref{fig:testsetup}(a). 
In abstract copies of our rigid body in Figs.~\ref{fig:testsetup}(b) and \ref{fig:frames}(b), they appear as blue cuboids.) 
Subfigures (a)--(d) in Fig.~\ref{fig:measurment} correspond, respectively, to the measurements from $A_1$--$A_4$. 
Each subfigure shows graphs of three (discrete) functions. 
For example, subfigure (d) shows  graphs of ${}^4\!\alpha_{i}\ag{\cdot}$, $i=1,2,3$. 
These functions are related to $A_4$'s acceleration vector in the following manner. 
Say  $^4\!\absf[\tau]$ is the non-dimensional acceleration of $A_4$. 
Then ${}^{4}\!\alpha_i[\tau]:={}^4\!\absf\ag{\tau}\cdot{}^{4}\!\ebsf_i\ag{\tau}$, $i=1,2,3$. 
The vectors  $\pr{^{4}\!\ebsf_i\ag{\tau}}_{i\in I}$  give the directions of accelerometer $A_4$'s measurement axes at the time instance $\tau$\footnote{
The dimensional counterparts of $\pr{^{4}\!\ebsf_i\ag{\cdot}}_{i\in I}$ at the initial and at a later time instance are, respectively, shown in Figs.~\ref{fig:testsetup}(b.i) and \ref{fig:testsetup}(b.ii). 
In each figure they appear as differently colored arrows emanating from $A_4$. 
Recall that the  accelerometers are rigidly fixed to the body and move with it. 
Therefore, the evolution of $\pr{^{4}\!\ebsf_i\ag{\tau}}_{i \in I}$ is the same as that in \eqref{eq:bfevolution} with the body basis vector $\ebsf_i\ag{\tau}$ replaced with $^{4}\!\ebsf_i\ag{\tau}$, and the laboratory basis vector $\ebsf_i$ replaced with the initial value of $^{4}\!\ebsf_i\ag{\cdot}$. 
In fact, the measurement axes of all four accelerometers, $\pr{^{\ell}\!\ebsf_i\ag{\tau}}_{i\in I}$, $\ell=1,2,3,4$, evolve similarly. 
Their evolution is illustrated in  Fig.~\ref{fig:testsetup}(b), and Fig.~\ref{fig:frames}(b). 
}. 
The functions graphed in the other subfigures (resp., $\pr{^{1}\!\alpha_i\ag{\cdot}}_{i\in I}$, $\pr{^{2}\!\alpha_i\ag{\cdot}}_{i\in I}$, and $\pr{^{3}\!\alpha_i\ag{\cdot}}_{i\in I}$) are similarly related to the acceleration vector evolutions at the other points (resp., $^{1}\!\absf\ag{\cdot}$, $^{2}\!\absf\ag{\cdot}$, and $^{3}\!\absf\ag{\cdot}$). 

\subsection{Applying the $\sqrt{\rm AO}$-algorithm}
\label{sec:ApplyingAO}
To recall, the $\sqrt{\rm{AO}}$-algorithm allows us to determine the angular velocity of a rigid body from acceleration measurements. 
The $\sqrt{\rm AO}$-algorithm takes as input discrete forms of the accelerometer measurements $^{\ell}\!\boldsymbol{\alpha}\ag{\cdot}:=\pr{^{\ell}\!\alpha_i\ag{\cdot}}_{i\in I}$, $\ell=1,2,3,4$, which we discussed in the previous section. 
Additionally, it requires the accelerometers' relative positions and orientations. 
These can be computed from the initial values of the accelerometers' position vectors\footnote{The dimensional form of $^{4}\xbsf\ag{\cdot}$ at the initial and a later time instance are shown in Fig.~\ref{fig:testsetup}(b.i) and (b.ii), respectively. 
} $^{\ell}\xbsf[0]$, $\ell=1,2,3,4$, and their measurement axes $\pr{^{\ell}\!\ebsf_i\ag{0}}_{i\in I}$. 
Our rigid body was 3D printed using fused deposition molding (FDM, material: polylactic acid (PLA), 3D printer: Original Prusa i3 MK3, Prusa Research, the Czech Republic).  
The parameters in its design were used to compute the initial values  $^{\ell}\xbsf[0]$, and  $\pr{^{\ell}\!\ebsf_i\ag{0}}_{i\in I}$. 
These values for our punctuated rotation trial are given in the caption of Fig.~\ref{fig:testsetup}. 

\begin{figure*}[t]
\centering
\includegraphics[width=15cm]{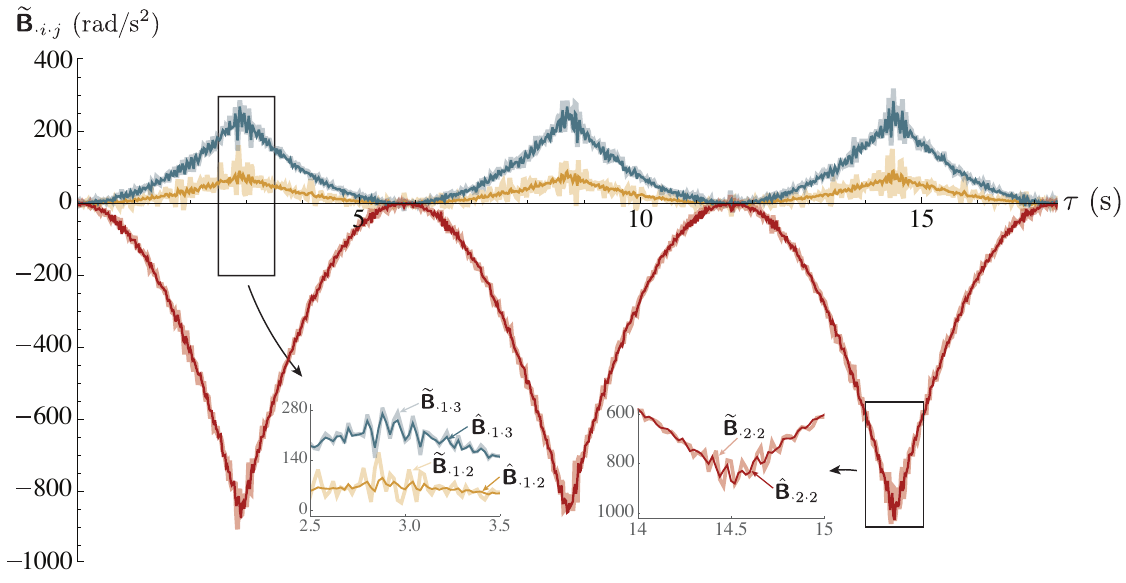}
    \caption{
	Graphs of certain components of $\Bbsf\ag{\cdot}$ in the punctuated rotation experiment. 
	To avoid clutter we arbitrarily chose to only show the $\pr{1,2}$, $\pr{1,3}$, and $\pr{2,2}$ components. 
	Each component is graphed twice: once using the raw measurement for  $\Bbsf\ag{\cdot}$, namely $\widetilde{\Bbsf}\ag{\cdot}$, which is computed directly from the acceleration measurements shown in Fig.~\ref{fig:measurment}, using a thick, solid, partially transparent line; and again using the best approximant of $\widetilde{\Bbsf}\ag{\cdot}$, namely $\widehat{\Bbsf}\ag{\cdot}$, which is obtained by applying Theorem~\ref{thm:main} to $\widetilde{\Bbsf}\ag{\cdot}$, using a thin, solid, opaque line. 
See \S\ref{sec:ApplyingAO} for details.
  } 
 \label{fig:Bmatrix}
\end{figure*}

On feeding  in the accelerometer measurements $^{\ell}\!\boldsymbol{\alpha}\ag{\tau_i}$, $i=1,2,\ldots$, and other experimental parameters, the $\sqrt{\rm AO}$-algorithm yields a discrete version of the function $\mathbb{R} \ni \tau \mapsto \Pbsf\ag{\tau}\in \Mvspace{3}$ (for details see \cite[Eq. (4.1)]{wan2022determining}). 
The body angular velocity matrix of the rigid body at the time instance $\tau$, $\overline{\Wbsf}\ag{\tau}$, is related to this output as
\begin{equation}
r\ag{\Bbsf\ag{\tau}}\ag{\overline{\Wbsf}\ag{\tau}}=\u{0},
\label{eq:ApplyAO1}
\end{equation}
where $\Bbsf\ag{\tau}  \in\Ktwoset{3}$ is the symmetric part of $\Pbsf\ag{\tau}$,  and $r\ag{\cdot}\ag{\cdot}$ is defined in \eqref{eq:symmp}. 
As discussed in \S\ref{sec:MotivationApplication}, due to experimental noise and errors, we do not obtain  $\Bbsf\ag{\tau}$ directly, but only an approximation for it, $\widetilde{\Bbsf}\ag{\tau}$, which usually does not lie in $\Ktwoset{3}$. 
(More specifically, we do not obtain $\Pbsf\ag{\tau}$ but only an approximation for it, $\widetilde{\Pbsf}\ag{\tau}$, and take the
$\text{symmetric part of} ~\widetilde{\Pbsf}\ag{\tau}=:\widetilde{\Bbsf}\ag{\tau}$ as an approximation for  $\Bbsf\ag{\tau}$.
Graphs of some of the component functions of $\widetilde{\Bbsf}\ag{\cdot}$  from the trial of our punctuated rotation experiment are shown in Fig.~\ref{fig:Bmatrix}.) 
The next step of the $\sqrt{\rm AO}$-algorithm involves finding $\widehat{\Bbsf}\ag{\tau_i}$, a best approximant for $\widetilde{\Bbsf}\ag{\tau_i}$ in $\Ktwoset{3}$. 
We calculate $\widehat{\Bbsf}\ag{\tau_i}$  by applying Algorithm~\ref{algo:ustar}. 
Specifically, we compute $\widehat{\Bbsf}\ag{\tau_i}$ as $\Ubsf^{\star}\ag{\widetilde{\Bbsf}\ag{\tau_i}}$, which in the current case can be written as
\begin{subequations}
\begin{align}
	&\Ubsf^{\star}\ag{\widetilde{\Bbsf}\ag{\tau_i}}
		=
			\Nbsf\ag{\tau_i} \DIAG_3\ag{0,\mu_1^{\star}\ag{\tau_i},\mu_1^{\star}\ag{\tau_i}} \Nbsf\ag{\tau_i}^\Tsf,
	\label{equ:Bdecomposition}
	\intertext{where}
	&\mu^{\star}_1\ag{\tau_i} = \left\{\begin{array}{l l}
	\pr{\lambda_{2}\ag{\tau_i} + \lambda_{3}\ag{\tau_i}}/2, & \quad \lambda_{2}\ag{\tau_i} + \lambda_{3}\ag{\tau_i} \leq 0, \\
	0, & \quad \text{otherwise},
	\end{array} \right.
\end{align}
\label{eq:Bhat3}
and $\lambda_2\ag{\tau_i}$, and $\lambda_3\ag{\tau_i}$ are the two smallest eigenvalues of $\widetilde{\Bbsf}\ag{\tau_i}$. 
The matrix $ \Nbsf\ag{\tau_i} \in \Mvspace{3} $ is constructed by taking the eigenvectors of $ \widetilde{\Bbsf}\ag{\tau_i} $ to be its columns, with the second and third columns corresponding to $ \lambda_2\ag{\tau_i} $ and $ \lambda_3\ag{\tau_i} $, respectively.
\end{subequations}
The raw matrices $\widetilde{\Bbsf}\ag{\tau_i}$ and their respective approximation  $\widehat{\Bbsf}\ag{\tau_i}$ are compared in  Fig.~\ref{fig:Bmatrix}. 

The final step of the $\sqrt{\rm AO}$-algorithm involves computing the approximate body angular velocity matrix $\widetilde{\overline{\Wbsf}}\ag{\tau_i}$ from the roots of $r\ag{\widehat{\Bbsf}\ag{\tau_i}}\ag{\cdot}$. 
For the current case the roots of $r\ag{\widehat{\Bbsf}\ag{\tau_i}}\ag{\cdot}$ come out to be
\begin{equation}
\widetilde{\overline{\Wbsf}}\ag{\tau_i}=\pm \Nbsf\ag{\tau_i}
\begin{pmatrix} 0 & 0 & 0 \\ 0 & 0 & -\sqrt{-\mu_1^{\star}\ag{\tau_i}} \\ 0 & \sqrt{-\mu_1^{\star}\ag{\tau_i}} & 0 \end{pmatrix}
\Nbsf\ag{\tau_i}^\Tsf
\label{equ:estimatedW}
\end{equation}
(see \cite[\S4.2]{wan2022determining} for details). 
The  positive and negative signs in \eqref{equ:estimatedW} correspond to the two roots of $r\ag{\widehat{\Bbsf}\ag{\tau_i}}\ag{\cdot}$. 
We choose the root that is closer to $\widetilde{\overline{\Wbsf}}\ag{\cdot}$'s value from  the previous time step (see \cite[\S4.2]{wan2022determining} for details).
Following \eqref{eq:wbf}, we compute an approximate body angular velocity vector at the time instance $\tau_i$--- $\widetilde{\busf{w}}\ag{\tau_i}$---as
$$
\widetilde{\busf{w}}\ag{\tau_i}=\star\widetilde{\overline{\Wbsf}}\ag{\tau_i}.
$$

Graphs of the component functions of $\widetilde{\busf{w}}\ag{\cdot}$ in the punctuated rotation trial are shown in Fig.~\ref{fig:motorcom}.

\subsection{Comparison}
\label{sec:comparison}
\begin{figure*}[t!]
\centering
\includegraphics[height=19.5cm]{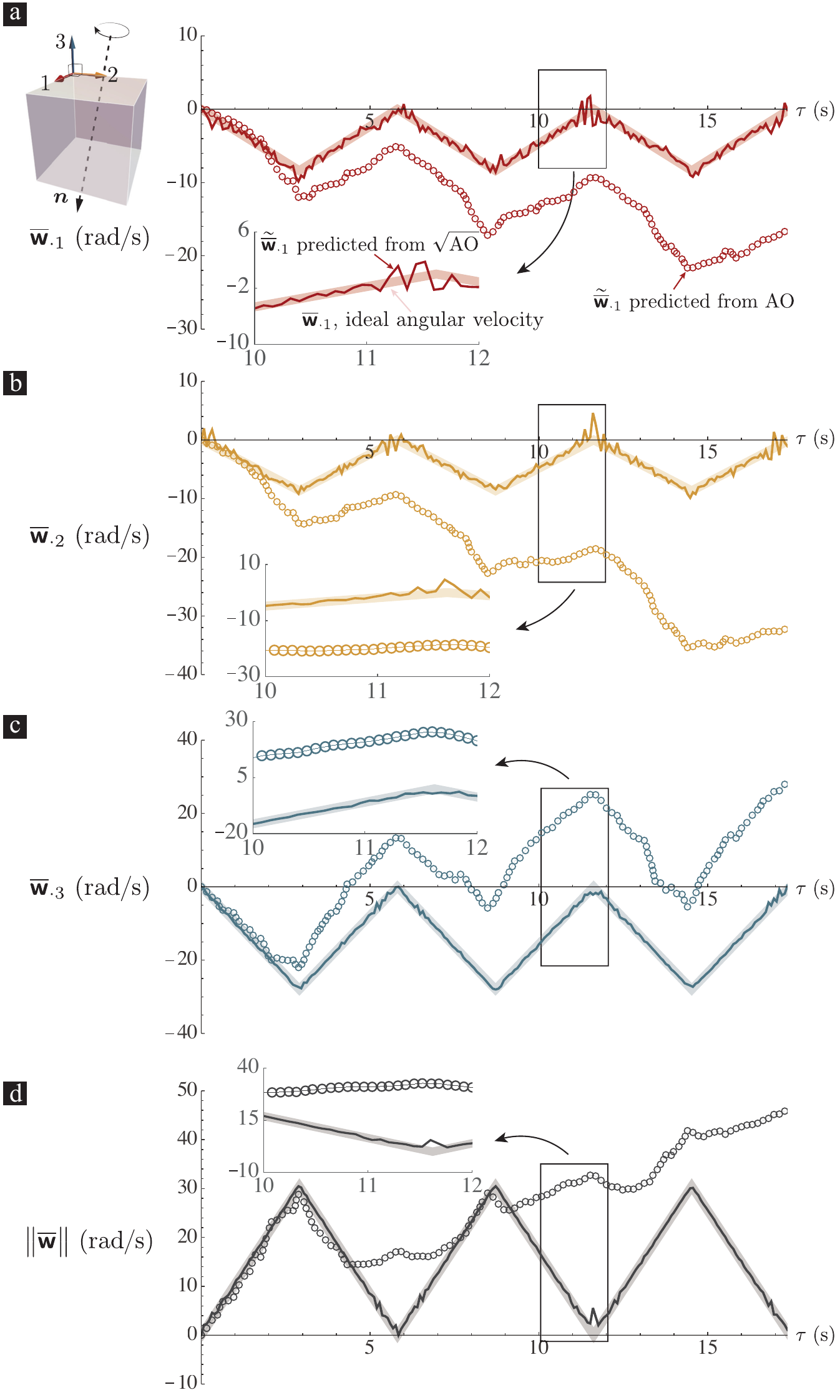}
    \caption{
	Body angular velocity in the punctuated rotation experiment. 
	Subfigures (a)--(c), respectively, show graphs of the body angular velocity's three components. (These are equivalently the components of the angular velocity w.r.t.~the body basis vectors, which are attached to the rigid body and rotate with it.
	The body basis vectors are illustrated, e.g., by the three arrows attached to the top corner of the gray cuboid (the rigid body) in the inset in (a).) 
	Each component is graphed thrice: once from the ideal body angular velocity, specifically the $\busf{w}\ag{\cdot}$ given by \eqref{eq:wbfIdeal} and \eqref{eq:PunctuatedRotationTheta} for $\omega_m = 31.41$ and $\tau_1 = 5.81$, using a thick, solid, slightly transparent line; once from the body angular velocity yielded by the $\sqrt{\rm AO}$ algorithm, namely $\widetilde{\busf{w}}\ag{\cdot}$, using  a thin, solid, opaque line; and finally from the body angular velocity yielded by the AO-algorithm, namely AO-$\widetilde{\busf{w}}\ag{\cdot}$, using hollow markers. 
	See  \S\ref{sec:comparison} for details. 
	Subfigure (d) shows the magnitude of the angular velocity. 
	It too is similarly graphed thrice. 
	}
 \label{fig:motorcom}
\end{figure*}

In addition to the graphs of $\widetilde{\busf{w}}\ag{\cdot}$'s components, Fig.~\ref{fig:motorcom} also shows the graphs of the ideal body angular velocity's components, i.e.~the components of ${\busf{w}}\ag{\cdot}$ given in \eqref{eq:wbfIdeal}, where $\theta\ag{\cdot}$ is given by \eqref{eq:PunctuatedRotationTheta} with $\omega_m = 31.41$ and $\tau_1 = 5.81$. 
The ideal angular velocity is what we expect to measure under ideal conditions, i.e., perfect rigid body motion, knowledge of experimental parameters, and acceleration measurements.  
As shown in Fig.~\ref{fig:motorcom}, the angular velocity components measured using the $\sqrt{\rm AO}$-algorithm closely align with their ideal counterparts. 
However, they exhibit some noise, particularly at time instances where the ideal components are non-smooth ($\tau=2.91$, 5.81, 8.72, etc. in Fig.~\ref{fig:motorcom}). 
Potential sources for this noise include accelerometer measurement noise, as well as spurious dynamics introduced by artifacts in the motor’s control system, chassis vibrations, and dynamic elastic deformations between the accelerometer particles.

We quantify the discrepancy between the measured and ideal angular velocities using the error metric
\begin{equation}
    \frac{\norm{\widetilde{\busf{w}}\ag{\cdot}-\busf{w}\ag{\cdot}}_{L^2(0,T)}}{\norm{\busf{w}\ag{\cdot}}_{L^2(0,T)}}.
	\label{eq:ErrorMetric}
\end{equation}
The  norm $\lVert \cdot \rVert_{L^2(0,T)}$, where $T>0$, is the standard $L^2$ function norm:
\begin{equation}
\norm{ \busf{w}\ag{\cdot} }_{L^2(0,T)}=\pr{\int_{0}^{T} \lVert \busf{w}\ag{\tau}\rVert^2}^{1/2},
\label{eq:ErrorMetric1}
\end{equation}
 and $[0,T]$  is $\busf{w}\ag{\cdot}$'s domain. 
Inside the integral in \eqref{eq:ErrorMetric1}, $\lVert \cdot \rVert$ denotes the standard Euclidean norm on $\mathbb{R}^3$. 
For the punctuated rotation experiment, the error is $0.062$ (or $6.2\%$). 
The error metrics for our oscillatory and constant rotation experiments are discussed in the \S\ref{sec:addrotationtest}. 

For reference, we also estimated the body angular velocity  using an alternative technique. 
In this technique  the body angular velocity is estimated by time integrating the experimental body angular acceleration. 
It can be shown that
\begin{equation}
	\busf{w}'\ag{\tau}=\star \pr{\text{skew part of the }{\usf{P}}\ag{\tau}},
	\label{eq:AOalgo}
\end{equation} where $\busf{w}'\ag{\cdot}$ is  the body angular acceleration and is the derivative of $\busf{w}\ag{\cdot}$, which, of course, is the body angular velocity. 
The matrix valued function $\usf{P}\ag{\cdot}$ in \eqref{eq:AOalgo} is the same one that we discussed in the context of the $\sqrt{\rm AO}$-algorithm's first step (see paragraph containing \eqref{eq:ApplyAO1}). 
It depends on the acceleration measurements, and the relative positions and orientations of the accelerometers. 
For a derivation of \eqref{eq:AOalgo} see \cite[\S 3.2]{rahaman2020accelerometer}. 
Following \eqref{eq:AOalgo}, we also estimated the body angular velocity in the experiment by integrating the skew part of $\widetilde{\Pbsf}\ag{\cdot}$, which was obtained by applying the first step of the $\sqrt{\rm AO}$-algorithm to that experiment. 
We call this alternative estimate of the body angular velocity the  $\widetilde{\busf{w}}\ag{\cdot}$ \textit{predicted from {\rm AO}-algorithm}, or AO-$\widetilde{\busf{w}}\ag{\cdot}$ for short. 
The components of AO-$\widetilde{\busf{w}}\ag{\cdot}$ in the punctuated rotation experiment are graphed in Fig.~\ref{fig:motorcom}. 
As observed, the discrepancies between the components of AO-$\widetilde{\busf{w}}\ag{\cdot}$ and their respective ideal components are significantly larger than those of the body angular velocity components estimated using the $\sqrt{\rm AO}$-algorithm. 
Notably, the error metric for AO-$\widetilde{\busf{w}}\ag{\cdot}$ is 1.64 (or 164\%), which is about 25 times greater than that of the body angular velocity estimated using the $\sqrt{\rm AO}$-algorithm. 

\section{Concluding remarks}
\label{sec:ConcludingRemarks}

In this work, we have proved a result concerning the best approximation of square matrices within the set of matrices that are the square of a skew-symmetric matrix. 
An immediate application of our result is  in numerical methods that aim to reproduce rigid body motion from acceleration measurements. 

\begin{enumerate}

    \item
    The present work uses the Frobenius norm to define a topology on the space of matrices. 
    In other works such as~\cite{higham1988computing}, best approximation in other norms are considered, such as the induced $\ell^2$-norm: for $\Absf \in \Mvspace{n}$,
    $$
    \| \Absf \|_2 := \sup_{\xbsf \in \Rbb^n} \frac{\| \Absf \xbsf \|}{\| \xbsf \|}.
    $$ 
    For rigid body motion, studying the best approximation problem in other norms may yield different error properties, though we emphasize that all norms are equivalent on finite-dimensional vector spaces. 
    Additionally, for other applications of the present problem, other norms may be more natural. 
    Nevertheless, additional work is required to extend the results presented here to other norms, and it is not apparent whether closed-form solutions may exist for these other norms.  
    \item
    In the context of the $\sqrt{\rm{AO}}$-algorithm, a natural question is how close is the estimated angular velocity matrix $\widetilde{\overline{\Wbsf}}$ to the actual angular velocity matrix $\overline{\Wbsf}$ at each time instance $\tau_i$. In dimensions $n = 2$ and $n = 3$, we can show that when $\fprod{\overline{\Wbsf},\widetilde{\overline{\Wbsf}}} > 0$
    $$
    \left\| \overline{\Wbsf} - \widetilde{\overline{\Wbsf}} \right\|_F^4 \leq C_n \left\| \Bbsf - \widehat{\Bbsf} \right\|_F^2,
    $$
    where $C_2 = 2$ and $C_3 = 8$. Additionally, because $\widehat{\Bbsf}$ is a best approximation of $\widetilde{\Bbsf}$, we have $\| \Bbsf - \widehat{\Bbsf} \|_F \leq 2 \| \Bbsf - \widetilde{\Bbsf}\|_F$ and so
    $$
    \left\| \overline{\Wbsf} - \widetilde{\overline{\Wbsf}} \right\|_F^4 \leq 4 C_n \left\| \Bbsf - \widetilde{\Bbsf} \right\|_F^2.
    $$
    Meanwhile, it may be shown via counterexample that no $C_n$ exists for $n > 3$. Proof of these results will be presented in a forthcoming publication.

    \item Since Theorem~\ref{thm:main} makes it possible to use the $\sqrt{\rm AO}$-algorithm on real data, we demonstrated this application by estimating the angular velocity of a rigid body using only acceleration measurements from four particles.
	Angular velocities measured via the $\sqrt{\rm AO}$-algorithm and through naive time integration of angular accelerations are compared in Fig.~\ref{fig:motorcom}. 
	What is especially notable in Fig.~\ref{fig:motorcom} is that for the case in which the angular velocity is measured via time integration the errors grow with time; while in the case in which they are measured via the $\sqrt{\rm AO}$-algorithm the errors remain bounded in time\footnote{There are mathematical justifications for both these observations. We will present these justifications in our forthcoming publication.}. 
	Based on these observations it is reasonable to conclude that the $\sqrt{\rm AO}$-algorithm is a significant improvement over the naive way of measuring angular velocity using only material particle accelerations.
	In fact, considering the advantages of measuring angular velocity from particle accelerations, i.e., using accelerometers, rather than by using gyroscopes\footnote{We mentioned these advantages briefly in \S\ref{sec:MotivationApplication}, but a thorough discussion is found in \S\ref{sec:table}.}, the $\sqrt{\rm{AO}}$-algorithm  is potentially a significant advancement in measuring angular velocity in general. 
	Since measurement of angular velocity is important in various fields---ranging from robotics to biomechanics---and, again, since Theorem~\ref{thm:main} enables the $\sqrt{\rm{AO}}$-algorithm's critical step, it follows that Theorem~\ref{thm:main} and the mathematical theory that we developed in support of it are likely to a have a significant technological impact.
\end{enumerate}

\section*{Acknowledgments}
he authors gratefully acknowledge support from the Panther Program, Tiger Program, and the Office of Naval Research (Dr. Timothy Bentley) under grants N000142112044 and N000142112054.

\section*{Declaration of Competing Interest}
The authors declare that they have no known competing financial interests or personal relationships that could have appeared to influence the work reported in this paper.

\section*{Author contributions}
Y.W.~carried out the research and experiments; B.E.G.P.~performed the mathematical analysis; H.K.~designed the research; all authors contributed to the writing of the manuscript.

\appendix
\section{Commercial off-the-shelf MEMS gyroscopes and accelerometers}
\label{sec:table}
A survey of commercially available off-the-shelf MEMS gyroscopes and accelerometers, along with their respective specifications, is presented in Table~\ref{tab:memssensor}.

\begin{longtblr}[
  caption = {Commercial off-the-shelf MEMS gyroscopes and accelerometers},
  label = {tab:memssensor},
  note{a} = {BMI serial product is an inertial measurement unit (IMU) which includes one triaxial accelerometer and one triaxial gyroscope.},
]{colspec = {llllc}}
\hline[0.8pt]
Company & Product & Bandwidth (Hz) & Size (length$\times$width$\times$height~mm$^3$) & Weight (g) \\
\hline
\SetCell[c=5]{c}{Gyroscopes} \\
\hline
STMicroelectronics \cite{STMicroelectronic} & A3G4250D & 110 & $4\times4\times1.1$ &--- \\
 & I3G4250D & 110 & $4\times4\times1.1$ &--- \\
   & L3G3250A & 140 & $3.5\times3\times1$ &--- \\
   & L3GD20 & 100 & $4\times4\times1$ &--- \\
   & L3G462A & 110 & $4\times4\times1.1$ &--- \\
   & LPR410AL & 110 & $4\times5\times1$ &--- \\
Analog Devices \cite{analogdevices} &ADXRS910 & 201 & $10.30\times10.42\times3.58$ &--- \\
 & ADXRS290 & 480 & $4.5\times5.8\times1.2$ &--- \\
   & ADXRS645 & 2000 & $15\times8\times2.85$ &--- \\
   & ADIS16137 & 400 & $35.6\times44\times13.8$ &--- \\
   & ADXRS642 & 2000 & $6.85\times6.85\times3.8$ &--- \\
   & ADXRS646 & 1000 & $6.85\times6.85\times3.8$ &--- \\
   & ADXRS453 & 77.5 & $9\times9\times4$ &--- \\
   & ADXRS620 & 2500 & $6.85\times6.85\times3.8$ &--- \\
   & ADXRS623 & 3000 & $6.85\times6.85\times3.8$ &--- \\
   & ADIS16060 & 1000 & $8.35\times8.2\times5.2$ &--- \\
     Tronics’s Microsystems SA \cite{Tronic} &G4300& 200 & $12\times12\times5.5$ &--- \\
   &G4050 & 30 & $12\times12\times5.5$ &--- \\
   &G3300 & 200 & $19.6\times11.5\times2.9$ &--- \\
    Bosch Sensortec\TblrNote{a} \cite{Bosch} &BMI323& 563 & $2.5\times3.0\times0.83$ &--- \\
   &BMI270 & 751 & $2.5\times3.0\times0.8$ &--- \\
   &BMI088 & 523 & $3\times4.5\times0.95$ &--- \\
   &BMI085 & 684 & $3\times4.5\times0.95$ &--- \\
    Silicon Sensing \cite{SliconSen} &CRM100& 160 & $5.7\times4.8\times1.2$ &--- \\
   &CRS43& 24 & $29\times29\times18.4$ &--- \\
   &CRH03& 100 & $47\times33.5\times25.4$ &--- \\
   &CMS390& 190 & $10.4\times6.7\times2.7$ &--- \\
    Gladiator Technologies \cite{Gladiator} &G300D& 600 & $25.4\times25.4\times15.2$ &--- \\
   &G200 & 200 & $8190$ (volume) &--- \\
   &G150Z & 200 & $18000$ (volume) &--- \\
    Safran \cite{Safran} &STIM210& 262 & $44.8\times38.6\times21.5$ &52 \\
   &STIM277H & 262 & $44.8\times38.6\times21.5$ &52 \\
   &STIM202 & 262 & $44.75\times38.6\times20$ &55 \\
\hline
\SetCell[c=5]{c}{Accelerometers} \\
\hline
STMicroelectronics \cite{STMicroelectronic} & AIS328DQ & 500 & $4\times4\times1.8$ &--- \\
& H3LIS331DL & 500 & $3\times3\times1$ &--- \\
& IIS2DLPC & 2500 & $2\times2\times0.7$ &--- \\
& IIS3DWB & 6300 & $2.5\times3\times0.86$ &--- \\
& LIS2DS12 & 6400 & $2\times2\times0.86$ &--- \\
& LIS2DUXS12 & 2500 & $2\times2\times0.74$ &--- \\
Analog Devices \cite{analogdevices} &ADXL1004 & 24000 & $5\times5\times1.8$ &--- \\
 & ADXL1001 & 11000 & $5\times5\times1.8$ &--- \\
   & ADXL382 & 8000 & $2.9\times2.8\times0.87$ &--- \\
   & ADXL203 & 2500 & $5\times5\times2$ &--- \\
   & ADXL354 & 1900 & $6\times6\times2.25$ &--- \\
   & ADXL316 & 1600 & $4\times4\times1.45$ &--- \\
   & ADXL355 & 1000 & $6\times6\times2.25$ &--- \\
     Tronics’s Microsystems SA \cite{Tronic} &A3050& 120 & $12\times12\times5.5$ &--- \\
   &A3140 & 300 & $12\times12\times5.5$ &--- \\
   &A3150 & 300 & $12\times12\times5.5$ &--- \\
    Bosch Sensortec\TblrNote{a} \cite{Bosch} &BMI323& 1677 & $2.5\times3.0\times0.83$ &--- \\
   &BMI270 & 684 & $2.5\times3.0\times0.8$ &--- \\
   &BMI088 & 523 & $3\times4.5\times0.95$ &--- \\
   &BMI085 & 684 & $3\times4.5\times0.95$ &--- \\
    Silicon Sensing \cite{SliconSen} &CAS211& 170 & $10.4\times6\times2.2$ &--- \\
   &CAS215& 170 & $10.4\times6\times2.2$ &--- \\
   &CAS291& 170 & $10.4\times6.7\times2.7$ &--- \\
   &CAS295& 170 & $10.4\times6.7\times2.7$ &--- \\
    Gladiator Technologies \cite{Gladiator} &GA50& 300 & $8.9\times8.9\times3.2$ &0.7 \\
   &A300D & 800 & $25.4\times25.4\times16.5$ &19.25 \\
   &A40 & 140 & $9832$ (volume) &--- \\
    Safran \cite{Safran} &MS1000& 200 & $8.9\times8.9\times3.23$ &1.5 \\
   &MS1000T & 200 & $8.9\times8.9\times3.23$ &1.5 \\
   &SI1000 & 550 & $8.9\times8.9\times3.23$ &1.5 \\
   &VS1000 & 1150 & $8.9\times8.9\times3.23$ &1.5 \\
\hline[0.8pt]
\end{longtblr}

\newpage
\section{Rotation test}
\label{sec:addrotationtest}
\begin{figure}[tp]
\centering
\includegraphics[width=0.9\textwidth]{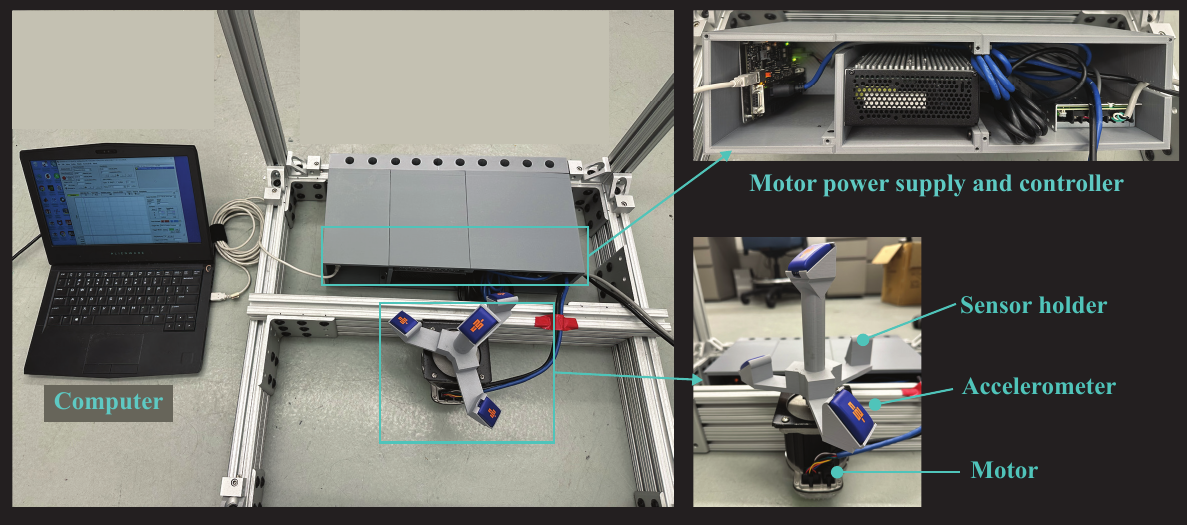}
\caption{Rotation test platform. Four triaxial accelerometers are mounted to a sensor holder. The sensor holder is rotating, driven by a feedback control BLDC motor, which is controlled by the computer.}
\label{fig:rotationtest}
\end{figure}

The setup of the rotation test platform is shown in Fig.~\ref{fig:rotationtest}.

\subsection{Additional rotation tests}

\begin{figure}[tp]
\centering
\includegraphics[width=14cm]{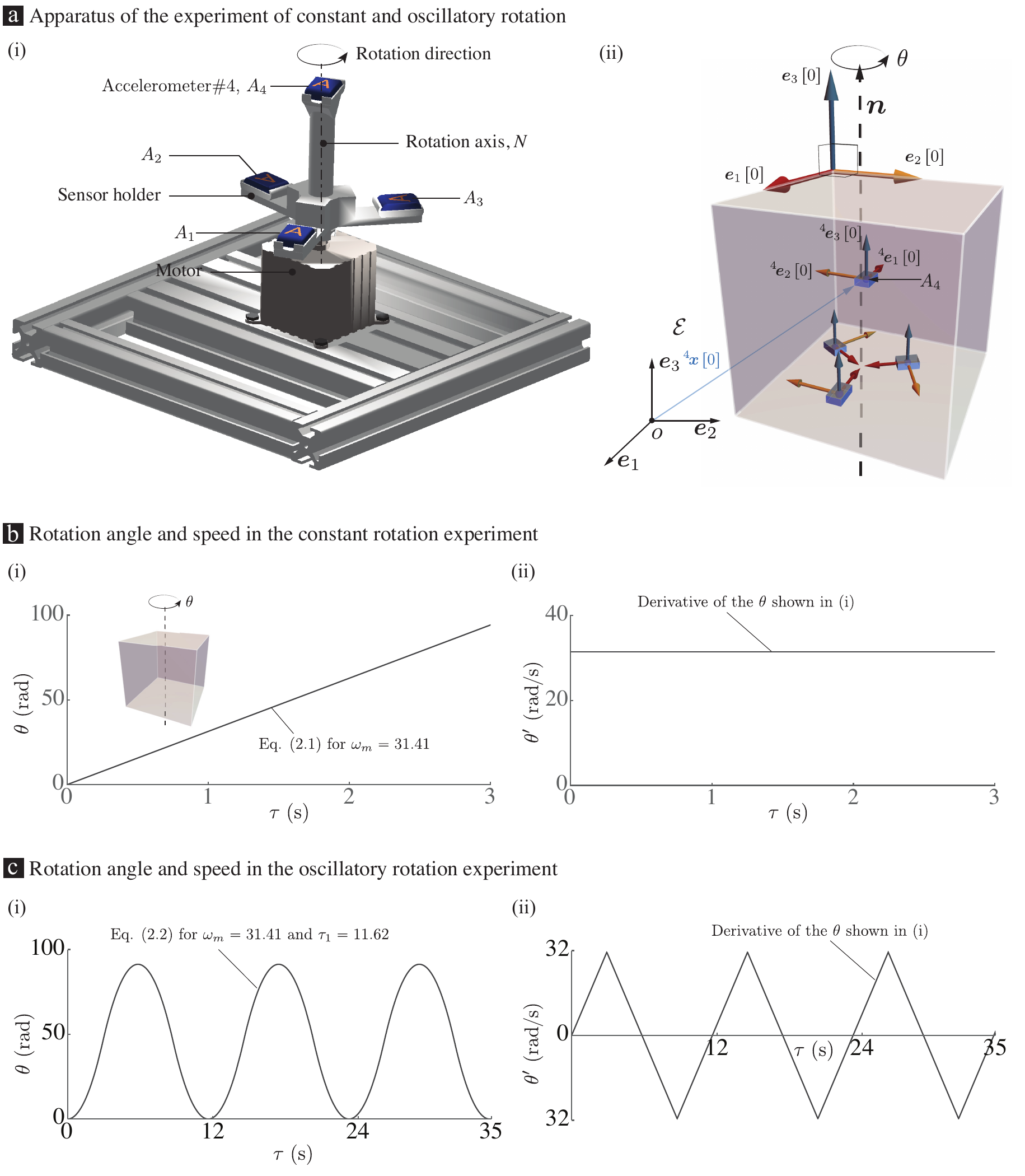}
    \caption{Set-up and geometry parameters of the constant and oscillatory rotation experiment. (a.i) shows a sketch of the rigid body rotation experiment. (a.ii) shows the initial configuration of the rigid body in both constant and oscillatory experiment.
    The components of four accelerometer's initial position vectors are ${}^1\!\usf{x}\ag{0}=\pr{0.064,0.001,-0.016}$,${}^2\!\usf{x}\ag{0}=\pr{-0.031,-0.055,-0.016}$, ${}^3\!\usf{x}\ag{0}=\pr{-0.032,0.055,-0.016}$, and ${}^4\!\usf{x}\ag{0}=\pr{-0.015,0.001,0.103}$. The components of $\pr{\lsc{1}\u{e}_i\ag{0}}_{i\in \pr{1,2,3}}$ are $\pr{\pr{-1,0,0},\pr{0,-1,0},\pr{0,0,1}}$; of $\pr{\lsc{2}\u{e}_i\ag{0}}_{i\in \pr{1,2,3}}$ are $\pr{\pr{0.50,0.86,0},\pr{-0.86,0.50,0},\pr{0,0,1}}$; of $\pr{\lsc{3}\u{e}_i\ag{0}}_{i\in \pr{1,2,3}}$ are $\pr{\pr{0.50,-0.86,0},\pr{0.86,0.50,0},\pr{0,0,1}}$; of $\pr{\lsc{4}\u{e}_i\ag{0}}_{i\in \pr{1,2,3}}$ are $\pr{\pr{-1,0,0},\pr{0,-1,0},\pr{0,0,1}}$; and of $\pr{\u{e}_i\ag{0}}_{i\in \pr{1,2,3}}$ are $\pr{\pr{1,0,0},\pr{0,1,0},\pr{0,0,1}}$.
    Subfigures (b.i) and (b.ii), respectively, shows the variation of the rotation angle and its derivative over time in the constant rotation experiment. Subfigures (c.i) and (c.ii), respectively, shows the variation of the rotation angle and its derivative over time in the oscillatory rotation experiment.
  }
\label{fig:testsetup_othertwotest}
\end{figure}

A sketch of the experimental setup for both constant and oscillatory rotation is shown in Fig.~\ref{fig:testsetup_othertwotest}(a.i). A drawing of a more abstract version of the rotation experiment is shown in Fig.~\ref{fig:testsetup_othertwotest}(a.ii).
The setups for the constant and oscillatory rotation experiments are the same, except for the different rotation profiles used.
In these two experiments, the rotation axis vector $\usf{n}$ was $\pr{0,0,1}$.
The initial values $\lsc{\ell}\usf{x}\ag{0}$, and $\pr{\lsc{\ell}\usf{e}_i\ag{0}}_{i\in I}$, $\ell=1,2,3,4$ are given in the caption of Fig.~\ref{fig:testsetup_othertwotest}.

\paragraph{Constant rotation}

\begin{figure*}[tp]
\centering
\includegraphics[height=19.5cm]{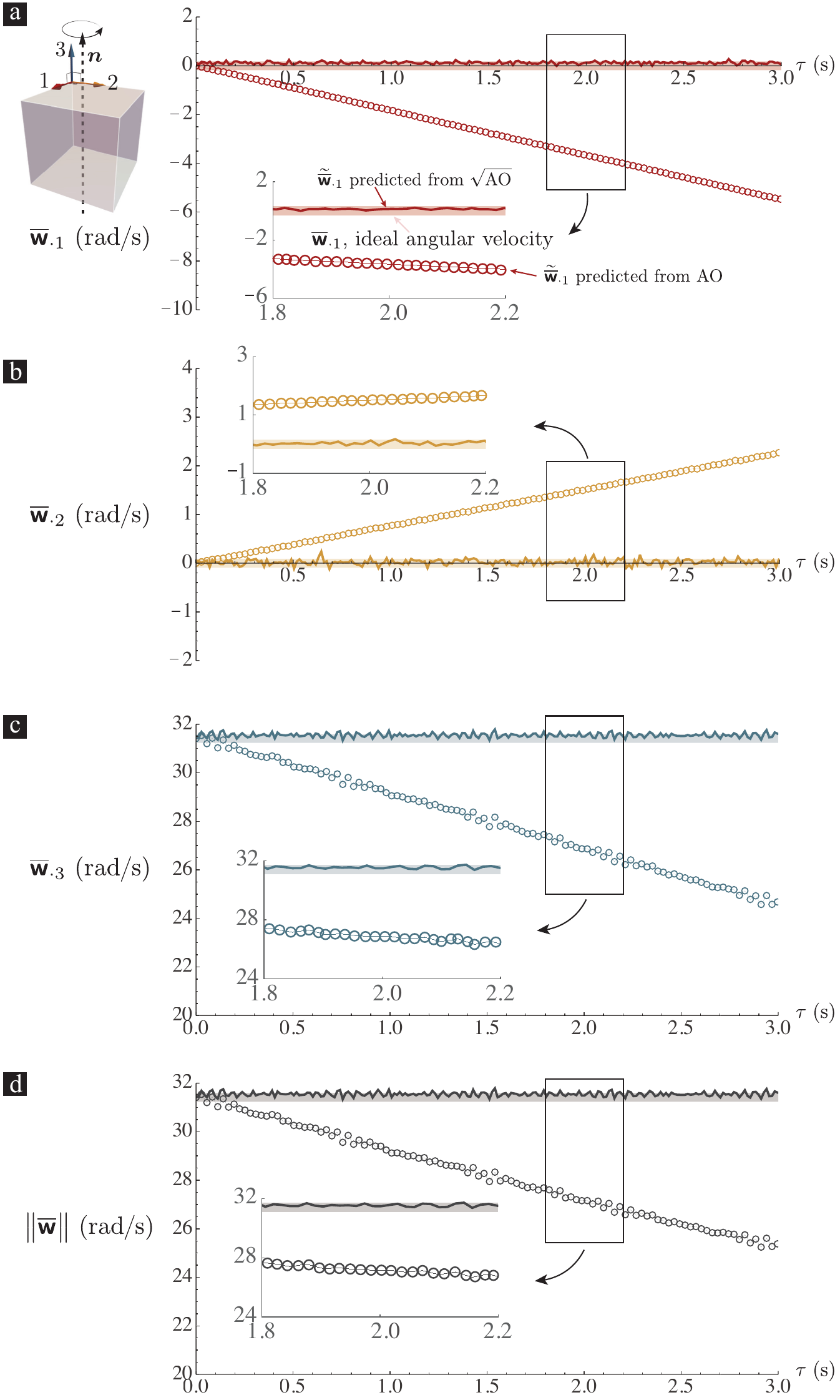}
    \caption{Body angular velocity in the constant rotation experiment. 
	Subfigures (a)--(c), respectively, show graphs of the body angular velocity's three components. (These are equivalently the components of the angular velocity w.r.t.~the body basis vectors, which are attached to the rigid body and rotate with it.
	The body basis vectors are illustrated, e.g., by the three arrows attached to the top corner of the gray cuboid (the rigid body) in the inset in (a).) 
	Each component is graphed thrice: once from the ideal body angular velocity, specifically the $\busf{w}\ag{\cdot}$ given by Eq. (32) in the main manuscript and \eqref{eq:theta_cons} for $\omega_m = 31.41$ using a thick, solid, slightly transparent line; once from the body angular velocity yielded by the $\sqrt{\rm AO}$ algorithm, namely $\widetilde{\busf{w}}\ag{\cdot}$, using  a thin, solid, opaque line; and finally from the body angular velocity yielded by the AO-algorithm, namely AO-$\widetilde{\busf{w}}\ag{\cdot}$, using hollow markers. 
	Subfigure (d) shows the magnitude of the angular velocity. 
	It too is similarly graphed thrice.
  }
 \label{fig:motorcom_constant}
\end{figure*}

In the constant rotation experiment, the rotation angle $\theta\ag{\cdot}$ is taken as
\begin{equation}
\theta\ag{\tau}=\omega_m\tau,
\label{eq:theta_cons}
\end{equation}
where $\omega_m=31.41$. Graphs of $\theta\ag{\cdot}$ in \eqref{eq:theta_cons} and its derivative are shown in Figs.~\ref{fig:testsetup_othertwotest}(b.i) and (b.ii), respectively.
The comparison of the body angular velocity components calculated using $\sqrt{\text{AO}}$-algorithm and AO-algorithm in the constant rotation experiment is shown in Fig.~\ref{fig:motorcom_constant}.
Fig.~\ref{fig:motorcom_constant} also shows the graphs of the ideal body angular velocity's components, i.e.~the components of ${\busf{w}}\ag{\cdot}$ given in Eq. (32) in the main manuscript, where $\theta\ag{\cdot}$ is given by \eqref{eq:theta_cons} with $\omega_m = 31.41$.

The error metric for the body angular velocity measured using $\sqrt{\text{AO}}$-algorithm in the trial of the constant rotation experiments is 0.0068 (or $0.68\%$), while the error metric for that using AO-algorithm  is 0.166 (or $16.60\%$).

\begin{figure*}[tp]
\centering
\includegraphics[height=19.5cm]{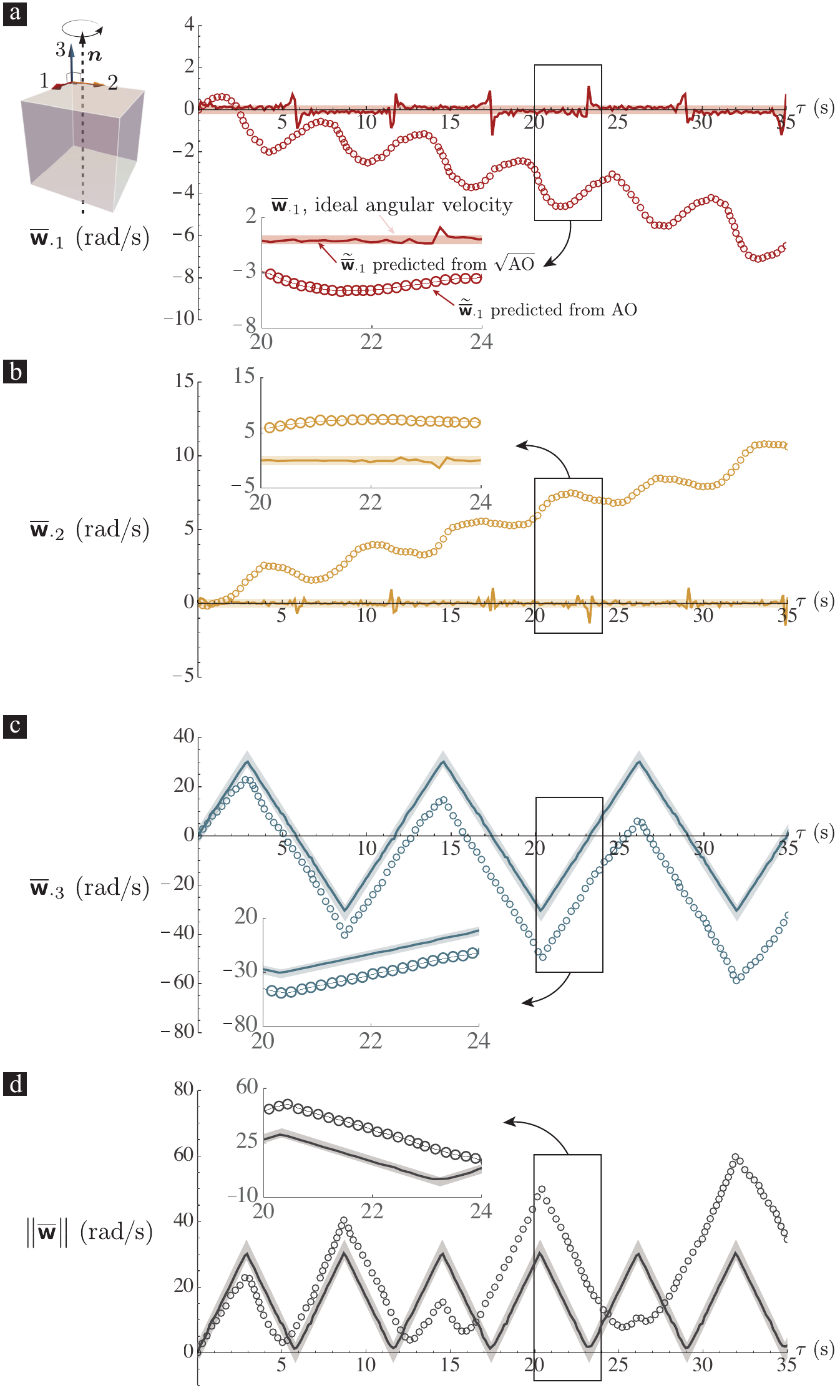}
    \caption{Body angular velocity in the oscillatory rotation experiment. 
	Subfigures (a)--(c), respectively, show graphs of the body angular velocity's three components. (These are equivalently the components of the angular velocity w.r.t.~the body basis vectors, which are attached to the rigid body and rotate with it.
	The body basis vectors are illustrated, e.g., by the three arrows attached to the top corner of the gray cuboid (the rigid body) in the inset in (a).) 
	Each component is graphed thrice: once from the ideal body angular velocity, specifically the $\busf{w}\ag{\cdot}$ given by Eq. (32) in the main manuscript and \eqref{eq:OscillatoryRotationTheta} for $\omega_m = 31.41$ and $\tau_1=11.62$ using a thick, solid, slightly transparent line; once from the body angular velocity yielded by the $\sqrt{\rm AO}$ algorithm, namely $\widetilde{\busf{w}}\ag{\cdot}$, using  a thin, solid, opaque line; and finally from the body angular velocity yielded by the AO-algorithm, namely AO-$\widetilde{\busf{w}}\ag{\cdot}$, using hollow markers. 
	Subfigure (d) shows the magnitude of the angular velocity. 
	It too is similarly graphed thrice.
  }
 \label{fig:motorcom_oscillatory}
\end{figure*}

\paragraph{Oscillatory rotation}
In the constant rotation experiment, the rotation angle $\theta\ag{\cdot}$ is taken as
\begin{equation}
    \theta[\tau]=  \frac{\omega_{\rm m}\tau_1}{8} +
    \frac{4\omega_{\rm m} \tau_1}{\pi^3} \sum_{n=1,3,5}
    \frac{(-1)^{(n+1)/2}}{n^3} \cos\ag{\frac{2n \pi \tau}{\tau_1}},
\label{eq:OscillatoryRotationTheta}
\end{equation}
where $\omega_m=31.41$ and $\tau_1=11.62$. Graphs of $\theta\ag{\cdot}$ in \eqref{eq:theta_cons} and its derivative are shown in Figs.~\ref{fig:testsetup_othertwotest}(c.i) and (c.ii), respectively.
The comparison of the component functions of body angular velocity calculated using $\sqrt{\text{AO}}$-algorithm and AO-algorithm are shown in Fig.~\ref{fig:motorcom_oscillatory}.
Fig.~\ref{fig:motorcom_oscillatory} also shows the graphs of the ideal body angular velocity's components, i.e.~the components of ${\busf{w}}\ag{\cdot}$ given in Eq. (32) in the main manuscript, where $\theta\ag{\cdot}$ is given by \eqref{eq:OscillatoryRotationTheta} with $\omega_m = 31.41$ and $\tau_1=11.62$.

The error metric for the body angular velocity measured via the $\sqrt{\text{AO}}$-algorithm in the trial of the oscillatory rotation experiments is 0.038 (or $3.80\%$), while via AO-algorithm is 1.0956 (or $109.56\%$).

\newpage
\section{Proofs of preliminary results}
\label{sec:proofs}
For completeness, we present proofs of the preliminary results below.

\subsection{Proof of Lemma~\ref{thm:frobenius-sym-skew}}

\ThmFrobeniusSymSkew*

\begin{proof}
\label{thm:frobenius-sym-skewproof}

\begin{align}
\| \Sbsf + \Kbsf \|_{F}^2 &= \langle \Sbsf + \Kbsf,\Sbsf + \Kbsf \rangle_F, \notag \\
&= \| \Sbsf \|_F^2 + \|\Kbsf\|_F^2 + 2 \langle \Sbsf,\Kbsf \rangle_F.
\label{eq:lemtst}
\end{align}
From the definition of the Frobenius inner product one can see that $\langle \Absf, \Bbsf \rangle_F = \langle \Absf^\Tsf, \Bbsf^\Tsf \rangle_F$. 
Hence,
$$
\begin{aligned}
\langle \Sbsf,\Kbsf \rangle_F &= \langle \Sbsf^\Tsf,\Kbsf^\Tsf \rangle_F,\\
&= \langle \Sbsf,-\Kbsf \rangle_F, \\
&=- \langle \Sbsf,\Kbsf \rangle_F,
\end{aligned}
$$
where the second equality is due to $\Kbsf$ being a skew-symmetric matrix. 
Therefore,  $\langle \Sbsf,\Kbsf \rangle_F  = 0$, and the result follows from \eqref{eq:lemtst}. 

\end{proof}

\subsection{Proof of Lemma~\ref{thm:frobenius-orthogonal}}

\ThmFrobeniusOrthogonal*

\begin{proof}
\label{thm:frobenius-orthogonalproof}
The Frobenius norm can be written in terms of the trace map $\TRACE\ag{\cdot}:\Mvspace{n}\to \mathbb{R}$, $$
\TRACE\ag{\Absf}=\sum_{i=1}^{n}\Absf_{\cdot i\cdot i},
$$
as
\begin{equation}
\|\Absf\|_F=\sqrt{\TRACE\ag{\Absf \Absf^{\Tsf}}}.
\end{equation}

To prove the first equality in the lemma, we express  $\| \Absf \Qbsf \|_F^2$ using the trace map as $\TRACE\ag{ \Absf \Qbsf (\Absf \Qbsf)^\Tsf }$. 
Which, owing to the property of the transpose map that $(\Absf \Qbsf)^\Tsf=\Qbsf^{\Tsf}\Absf^{\Tsf}$, simplifies to $\TRACE\ag{ \Absf \Qbsf \Qbsf^\Tsf \Absf^\Tsf}$. 
Since we know from Eq.~(4) of the main manuscript that $\Qbsf \Qbsf^\Tsf=\IdentityMatrix_n$, this expression further simplifies to $\TRACE\ag{\Absf \Absf^\Tsf}=\|$\Absf$\|_F^2$. 
The result follows on taking the positive square root of this equation. 

For showing the second equality in the lemma we follow a similar strategy and  write $\| \Qbsf \Absf \|_F^2$ in terms of the trace map as $\TRACE\ag{\Qbsf \Absf  (\Qbsf \Absf)^\Tsf }$, or equivalently as $\TRACE\ag{\Qbsf \Absf  \Absf^\Tsf \Qbsf^\Tsf}$. 
The trace map has the property that $\TRACE\ag{\Absf \Bbsf}=\TRACE\ag{\Bbsf \Absf}$,
where $\Absf, \Bbsf \in \Mvspace{n}$. 
Using this property and commuting the composition of the first three matrices, namely $\Qbsf \Absf  \Absf^\Tsf$, and the last matrix, namely $\Qbsf^\Tsf$, in $\TRACE\ag{\Qbsf \Absf  \Absf^\Tsf \Qbsf^\Tsf}$ we get  that $\| \Qbsf \Absf \|_F^2=\TRACE\ag{\Qbsf^\Tsf \Qbsf \Absf  \Absf^\Tsf }$. 
This equation simplifies to $\| \Qbsf \Absf \|_F^2=\|\Absf\|_F^2$, since  $\Qbsf^\Tsf \Qbsf=\IdentityMatrix_n$, and
$\TRACE\ag{\Absf  \Absf^\Tsf }=\|\Absf\|_F^2$. 
The result follows on taking the positive square root of this equation. 

\end{proof}

\subsection{Proof of Lemma~\ref{thm:K2-eigenvalues}}

\ThmKTwoEigenvalues*

\begin{proof}
\label{thm:K2-eigenvaluesproof}

Proof that $\Ktwoset{n} \subseteq \Vc_n$ may be found in~\cite[Appendix C]{wan2022determining}. 
The proof relies on the fact that the eigen-decomposition of $\Kbsf \in \Kvspace{n}$ has pure imaginary eigenvalues (which come in conjugate pairs) and orthogonal (complex) eigenvectors. 
Below we show that $\Vc_n \subseteq \Ktwoset{n}$. (This result can also be found in \cite{Rinehart1960}.) 
\begin{enumerate}
\item Let $\Absf \in \Vc_n$, i.e.,
$\Absf=\Nbsf \Dbsf \Nbsf^\Tsf$, where $\Nbsf \in \Oset{n}$, when $n$ is even (resp.~odd)
$\Dbsf=\DIAG_n\ag{\mu_1, \mu_1,\ldots, \mu_k, \mu_k}$ (resp.~$\Dbsf=\DIAG_n\ag{0,\mu_1, \mu_1,\ldots, \mu_k, \mu_k}$), and $k=n/2$ (resp.~$k=(n-1)/2$). 

\item When $n$ is even (resp.~odd) we  define
\begin{subequations}
\begin{align}\Kbsf = \DIAG\ag{\begin{pmatrix} 0 & -\sqrt{-\mu_1} \\ \sqrt{-\mu_1} & 0 \end{pmatrix}, \ldots, \begin{pmatrix} 0 & -\sqrt{-\mu_k} \\ \sqrt{-\mu_k} & 0 \end{pmatrix} }\footnotemark\\
\left(\text{resp. }\Kbsf = \DIAG\ag{ \begin{pmatrix} 0 \end{pmatrix} ,\begin{pmatrix} 0 & -\sqrt{-\mu_1} \\ \sqrt{-\mu_1} & 0 \end{pmatrix}, \ldots, \begin{pmatrix} 0 & -\sqrt{-\mu_k} \\ \sqrt{-\mu_k} & 0 \end{pmatrix} }\right).
\end{align}
\label{eq:Kdef}
\end{subequations}
\footnotetext{For square matrices $\Absf_1 \in \Mvspace{n_1},\ldots,\Absf_m \in \Mvspace{n_m}$ (which may have different sizes $n_1,\ldots,n_m$, respectively), we let $\DIAG\ag{\Absf_1,\ldots,\Absf_m}$ denote the block diagonal matrix
\begin{equation}
\begin{pmatrix}
\Absf_1 & \ZeroMatrix & \cdots & \ZeroMatrix \\
\ZeroMatrix & \Absf_2 & \cdots & \ZeroMatrix \\
\vdots & \vdots & \ddots & \vdots \\
\ZeroMatrix & \ZeroMatrix & \cdots & \Absf_m
\end{pmatrix} \in \Mvspace{\sum_i n_i},
\end{equation}
where the $\ZeroMatrix$ entries denote matrices whose entries are all zero (and whose dimensions we infer by context).} 
It is straightforward to show that $\Kbsf^\Tsf = - \Kbsf$, and hence $\Kbsf \in \Kvspace{n}$. 
Moreover, a direct calculation\footnote{A useful identity for carrying out this calculation is
\begin{equation}
\begin{pmatrix} 0 & -\mu \\ \mu & 0 \end{pmatrix}^2 = \begin{pmatrix} -\mu^2 & 0 \\ 0 & -\mu^2 \end{pmatrix},
\end{equation} where $\mu \in \mathbb{R}$.
} gives $\Kbsf^2 = \Dbsf$. 

\item Let $\usf{B}:= \Nbsf \Kbsf \Nbsf^\Tsf$, then
\begin{equation}
\usf{B}^\Tsf = (\Nbsf \Kbsf \Nbsf^\Tsf)^\Tsf = \Nbsf \Kbsf^\Tsf \Nbsf^\Tsf = -\Nbsf \Kbsf \Nbsf^\Tsf = - \usf{B}.
\label{eq:btranspose}
\end{equation}
The third equality in \eqref{eq:btranspose} follows from the fact that $\Kbsf$ is skew-symmetric. 
Furthermore,
\begin{equation}
\usf{B}^2 = (\Nbsf \Kbsf \Nbsf^\Tsf) (\Nbsf \Kbsf \Nbsf^\Tsf) = \Nbsf \Kbsf^2 \Nbsf^\Tsf = \Nbsf \Dbsf \Nbsf^{\Tsf} = \usf{A}.
\label{eq:bsquare}
\end{equation}
The second equality in \eqref{eq:bsquare} follows from the fact that $\Nbsf$ is orthogonal; the third from the fact that $\Kbsf^2=\Dbsf$ (see Step~2); and the last equality from $\Absf$'s definition. 
Since we know that $\usf{B} \in \Kvspace{n}$, from \eqref{eq:btranspose}, and $\usf{B}^2 = \usf{A}$, from \eqref{eq:bsquare}, it follows from  $\Ktwoset{n}$'s definition that $\usf{A} \in \Ktwoset{n}$. 
\end{enumerate}
\end{proof}

We briefly comment on the uniqueness of the matrix $\Bbsf$ in the proof of Lemma \ref{thm:K2-eigenvalues}. In particular, this matrix may not be unique. 
Taking the $\Kbsf$ given in \eqref{eq:Kdef}, multiplying one or more of its diagonal blocks with $-1$, and performing an orthogonal similarity transformation on it using $\Nbsf$ will yield a matrix which too lies in $\Kvspace{n}$ and squares to $\Absf$. 

\subsection{Proof of Lemma~\ref{thm:lambda-m-n}}

\ThmLambdaMN*

\begin{proof} Let $\Qbsf \in \Oset{n}$. The matrices $\Lambdabsf$ and $\Qbsf \Dbsf \Qbsf^\Tsf$ are both symmetric and have real eigenvalues $(\lambda_1,\ldots,\lambda_n)$ and $(\mu_1,\ldots,\mu_n)$ arranged in descending order. Applying~\cite[Corollary 6.3.8]{horn2012matrix} (where, using the notation from the reference, we take matrices $\Absf$ and $\Absf+\Ebsf$ to be $\Lambdabsf$ and $\Qbsf \Dbsf \Qbsf^\Tsf$, respectively) gives
$$
\sum_{i=1}^n (\lambda_i - \mu_i)^2 \leq \| \Lambdabsf - \Qbsf \Dbsf \Qbsf^\Tsf \|_F^2.
$$
The left-hand-side is identically $\| \Lambdabsf - \Dbsf \|_F^2$. Since this inequality holds for any $\Qbsf \in \Oset{n}$ (and is an equality in the particular case when $\Qbsf = \IdentityMatrix_n$), the proof follows.
\end{proof}

\subsection{Alternate Proof of Lemma~\ref{thm:lambda-m-n}}

We also present a second, longer proof which does not make use of~\cite[Corollary 6.3.8]{horn2012matrix}. Before stating the proof, we first discuss an isomorphism between the space of matrices and the space of $n$-tuples of vectors.

\paragraph{The Euclidean space $\Mprime{n}$}
Consider the set $\Mprime{n}:=\overbrace{\mathbb{R}^n \times \mathbb{R}^n\times \ldots \mathbb{R}^n}^{\text{$n$-times}}$. Let $\widetilde{\Xbsf} = \pr{\xbsf_1,\ldots,\xbsf_n}$ and $\widetilde{\Ybsf} = \pr{\ybsf_1,\ldots,\ybsf_n}$ belong to $\Mprime{n}$. Define the function $\Mpprod{\cdot,\cdot}: \Mprime{n} \times \Mprime{n} \to \mathbb{R}$,
\begin{equation}
\Mpprod{\widetilde{\Xbsf},\widetilde{\Ybsf}}=   \sum_{i=1}^{n}\langle \xbsf_{i},\ybsf_{i} \rangle.
\label{eq:newIP}
\end{equation}
It can be shown that $\Mprime{n}$ together with $\Mpprod{\cdot,\cdot}$ is a Euclidean vector space. In fact, the inner product spaces $\pr{ \Mvspace{n}, \fprod{\cdot, \cdot}}$ and $\pr{\Mprime{n}, \Mpprod{ \cdot, \cdot}}$ are isomorphic. For example, there exists the bijective, linear, inner product preserving map $\texttt{iso}\ag{\cdot}: \Mvspace{n}\to \Mprime{n}$,
\begin{equation}
    \texttt{iso}\ag{\Xbsf}=
    \pr{
    \pr{X_{j1}}_{j\in \Ic},
    \pr{X_{j2}}_{j\in \Ic},
    \ldots,
    \pr{X_{jn}}_{j\in \Ic}
    },
\end{equation}
where $X_{ij}$ is the $i\mbox{-}j^{\rm th}$ component of $\Xbsf$.
Note that the map $\texttt{iso}\ag{\cdot}$ is also a homeomorphism. The inverse of $\texttt{iso}\ag{\cdot}$ can be explicitly written as
\begin{equation}
   \texttt{iso}^{-1}\ag{\pr{\ybsf_1, \ldots, \ybsf_n}}=
\pr{\ybsf_i}_{i\in \Ic}^\Tsf.
\end{equation}


\begin{proof}\hfill
\begin{enumerate}
\item We begin by showing the minimization problem in \eqref{eq:minp} is equivalent to maximizing $g\ag{\cdot} : \Oset{n} \rightarrow \Rbb$, where
$$
g\ag{\Ybsf} = \fprod{ \Lambdabsf, \Ybsf \Dbsf \Ybsf^\Tsf }.
$$
Note that the function $g\ag{\cdot}$ is continuous, and that its domain, namely $\Oset{n}$, is a compact subset\footnote{For any $\Qbsf \in \Oset{n}$, we have $\| \Qbsf \|_F = \sqrt{n}$, which follows from Eq.~(1.2) of the main manuscript.
~Therefore $\Oset{n}$ is bounded. 
The set $\Oset{n}$ is also closed since it is the pre-image of the closed set $\left \{\ZeroMatrix\right\}\subset \Mvspace{n}$ under the continuous map $\Mvspace{n}\ni \Ybsf \mapsto \pr{\Ybsf^{\Tsf}\Ybsf-\IdentityMatrix_n}\in \Mvspace{n}$, where $\ZeroMatrix$ is the matrix whose entries are all zero. 
} of the Euclidean vector space $\Mvspace{n}$. 
Therefore, there exists a $\Ybsf^{\star}\in \Oset{n}$  such that
$$
\max_{\Ybsf\in \Oset{n}}g\ag{\Ybsf}=g\ag{\Ybsf^{\star}}.
$$

\begin{enumerate}
\item Via matrix algebra and the definition of the Frobenius inner product and norm,
\begin{align}
\| \Lambdabsf - \Qbsf \Dbsf \Qbsf^\Tsf \|_F^2 &= \fprod{ \Lambdabsf - \Qbsf \Dbsf \Qbsf^\Tsf, \Lambdabsf - \Qbsf \Dbsf \Qbsf^\Tsf }, \notag \\
&=\| \Lambdabsf \|_F^2 + \| \Qbsf \Dbsf \Qbsf^\Tsf \|_F^2 - 2 \fprod{ \Lambdabsf, \Qbsf \Dbsf \Qbsf^\Tsf }, \notag \\
&= \| \Lambdabsf \|_F^2 + \| \Dbsf \|_F^2 - 2 \fprod{ \Lambdabsf, \Qbsf \Dbsf \Qbsf^\Tsf }, \notag \\
&= \| \Lambdabsf \|_F^2 + \| \Dbsf \|_F^2 - 2 g\ag{\Qbsf}.
\label{eq:IPExpansion}
\end{align}
For the third equality, we used Lemma~\ref{thm:frobenius-orthogonal} to equate $ \| \Qbsf \Dbsf \Qbsf^\Tsf \|_F$ with $\| \Dbsf \|_F$. 

\item It follows from \eqref{eq:IPExpansion} and the definition of $\Ybsf^\star$ that
\begin{align}
   \min_{\Qbsf \in \Oset{n}} \| \Lambdabsf - \Qbsf \Dbsf \Qbsf^\Tsf \|_F=\| \Lambdabsf - \Ybsf^{\star} \Dbsf \left. \Ybsf^{\star}\right.^\Tsf \|_F.
\end{align} 



\end{enumerate}

\item\label{eq:NewGoal2} Completing the proof now requires us to show that $g\ag{\IdentityMatrix_n} = g\ag{\Ybsf^\star}$, i.e. $\IdentityMatrix_n$ is a maximizer of $g\ag{\cdot}$. For ease of computations, we will instead show that $\tilde{g}\ag{\cdot} = g\circ\texttt{iso}^{-1}\ag{\cdot}$, which is defined over the set
$$
\Oprime{n}=\left\{ \pr{\ybsf_1,\ldots,\ybsf_n}\in \Mprime{n}\SUCHTHAT \langle \ybsf_i,\ybsf_j\rangle=\delta_{ij},\ i,j\in \Ic \right\},
$$
achieves its maximum value at $\widetilde{\IdentityMatrix}_n = \texttt{iso}[\IdentityMatrix_n]$. By the use of isomorphism, $\IdentityMatrix_n$ is also a maximizer of $g\ag{\cdot}$.

\begin{enumerate}

\item A more explicit form of $\tilde{g}\ag{\cdot}$ is
\begin{align}
    \tilde{g}\ag{\pr{\ybsf_1,\ldots, \ybsf_n}}&=\langle \Lambdabsf, \texttt{iso}^{-1}\ag{\pr{\ybsf_1, \ybsf_2, \ldots, \ybsf_n}}\, \Dbsf\, \texttt{iso}^{-1}\ag{\pr{\ybsf_1, \ybsf_2, \ldots, \ybsf_n}}^\Tsf \rangle_F,\notag\\
    &= \sum_{i=1}^n \mu_i \langle \ybsf_i, \Lambdabsf \ybsf_i\rangle,
\label{eq:f}
\end{align}
where we recall that each $\ybsf_i\in \mathbb{R}^n$.  

\end{enumerate}

\item \label{thm:lambda-m-n:step1} We first show Step \ref{eq:NewGoal2} for the case where $\lambda_1 > \lambda_2 > \ldots > \lambda_n$ and $\mu_1 > \mu_2 > \ldots > \mu_n$ and $\lambda_i,\mu_i \neq 0$ for all $i \in \Ic$. 
\begin{enumerate}
\item Let $\Ytbsf^{\star}$ be a global maximum of $\tilde{g}\ag{\cdot}$. The set $\Oprime{n}$ is a smooth manifold\forus{Since $\Oprime{n}$ is smooth manifold we can talk about derivative and gradients on it. This enables us to now apply Fermat's theorem to this problem $\tilde{\mathcal{P}}$. Firstly we note that since $\Oprime{n}$ is a manifold it does not have a boundary. } 
and $\tilde{g}\ag{\cdot}$ is smooth everywhere on $\Oprime{n}$.\forus{Since $\tilde{g}\ag{\cdot}$ is differentiable everywhere on $\Oprime{n}$ it is clearly differentiable at points on $\Oprime{n}$ that are not the boundary points.} 
Hence, $\tilde{\Ybsf}^{\star}$ is a stationary point of $\tilde{g}\ag{\cdot}$. 

\item \forus{Without proof or justification we state that}A necessary condition for $\pr{\ybsf^{\star}_1,\ybsf^{\star}_2,\ldots, \ybsf^{\star}_n}:=\Ytbsf^{\star}$ to be a stationary point of $\tilde{g}\ag{\cdot}$ is that for each $i \in \Ic$
\begin{subequations}
\label{eqs:StationarityConditions}
\begin{equation}
\partial_i \tilde{g}\ag{\Ytbsf^{\star}} +\sum_{m=1}^{n}\sum_{k=1}^{n}L_{mk}\partial_i H_{m k}\ag{\Ytbsf^{\star}} = \ZeroMatrix\in \mathbb{R}^n,
\label{eq:stationarycond1}
\end{equation}
where
\begin{equation}
H_{ij}\ag{\ybsf_1, \ybsf_2, \ldots, \ybsf_n}:=
\langle
\ybsf_i, \ybsf_j\rangle - \delta_{ij},
\label{eq:H}
\end{equation}
$i, j\in \Ic$. 
In \eqref{eq:stationarycond1}, the quantities $\partial_i\tilde{g}\ag{\Ytbsf^{\star}}$, $\partial_i H_{m k}\ag{\Ytbsf^{\star}}$, respectively, denote
\begin{align}
    \pder{\tilde{g}\ag{\Ytbsf^{\star}}}{\ybsf_i}&:=\pr{\pder{\tilde{g}\ag{\Ytbsf^{\star}}}{\ytbsf_{i\cdot 1}},\ldots, \pder{\tilde{g}\ag{\Ytbsf^{\star}}}{\ \ybsf_{i\cdot n}}},\\
    \pder{H_{m k}\ag{\Ytbsf^{\star}}}{\ybsf_i}&:=\pr{\pder{ H_{mk}\ag{\Ytbsf^{\star}}}{ \ybsf_{i\cdot 1}},\ldots, \pder{ H_{mk}\ag{\Ytbsf^{\star}}}{ \ybsf_{i\cdot n}}}.
\end{align}
\end{subequations} 
The quantities $L_{mk}$, $m, k\in \Ic$, in \eqref{eq:stationarycond1} are components of a real symmetric matrix. 
That is,  $L_{mk}=\Lbsf_{\cdot m\cdot k}$, where $\Lbsf \in \Svspace{n}$. 

\item  It follows from \eqref{eq:f} and \eqref{eqs:StationarityConditions} that
\begin{equation}
 2 \mu_i \Lambdabsf \ybsf^{\star}_i - 2 \sum_{k=1 }^n L_{ik} \ybsf^{\star}_k = \bm{\mathsf{0}},
\label{eq:stationarycond}
\end{equation}
holds for each $i\in \Ic$. 
Taking the inner product of \eqref{eq:stationarycond} with $\ybsf^{\star}_j$,  noting that
$\langle {\ybsf^{\star}_k}, \ybsf^{\star}_j \rangle= \delta_{kj}$, and simplifying we get that
\begin{equation}
L_{ij} = \mu_i \langle \ybsf^{\star}_j, \Lambdabsf \ybsf^{\star}_i\rangle,
\label{eq:Lij}
\end{equation}
where $i, j\in \Ic$. 
Computing $L_{ij}$ and $L_{ji}$ from \eqref{eq:Lij} and setting them equal to each other, since they are the components of a symmetric matrix,
we get that
\begin{align}
    \mu_i \langle \ybsf^{\star}_j, \Lambdabsf \ybsf^{\star}_i\rangle &=\mu_j \langle \ybsf^{\star}_i, \Lambdabsf \ybsf^{\star}_j\rangle,\notag \\
    &=\mu_j \langle \Lambdabsf^{\Tsf} \ybsf^{\star}_i,  \ybsf^{\star}_j\rangle,\notag \\
    &=\mu_j \langle \Lambdabsf \ybsf^{\star}_i,  \ybsf^{\star}_j\rangle.
    \label{eq:Lij1}
\end{align} 
The second equality in \eqref{eq:Lij1} follows from the definition of the matrix transpose. 
The third equality in \eqref{eq:Lij1} follows because $\Lambdabsf$ is a diagonal matrix, and hence a symmetric matrix. 
Thus, by the symmetry of the inner products, $\pr{\mu_i -\mu_j}\langle \ybsf^{\star}_j, \Lambdabsf \ybsf^{\star}_i\rangle =0$ for any $i, j\in \Ic$. 
By assumption, $\mu_i \neq \mu_j$ when $i\neq j$, and so
\begin{equation}
 \langle \ybsf^{\star}_j, \Lambdabsf \ybsf^{\star}_i\rangle = 0,
 \label{eq:Lij3}
\end{equation}
for $i \neq j$. 
Along with \eqref{eq:Lij}, we must have that $L_{ij}=0$ for $i\neq j$. 

\item Given $L_{ij}=0$ for $i\neq j$ and $L_{ii} = \mu_i \langle\ybsf_i^\star,\Lambdabsf \ybsf_i^\star \rangle$, \eqref{eq:stationarycond} may be simplified to
\begin{equation}
 2 \mu_i \left( \IdentityMatrix_n - \ybsf^{\star}_i\otimes {\ybsf^{\star}_i} \right) \Lambdabsf \ybsf^{\star}_i = \ZeroMatrix,
\label{eq:stationarycond3}
\end{equation}
where $\ybsf_i^{\star}\otimes \ybsf_i^{\star}\in \Mvspace{n}$ such that $\pr{ \ybsf_i^{\star}\otimes \ybsf_i^{\star}}_{\cdot m\cdot n }=\ybsf_{i\cdot m}^{\star}\ybsf_{i\cdot n}^{\star}$. 
Because $\mu_i \neq 0$, 
it must either be the case that $\Lambdabsf \ybsf^{\star}_i = \ZeroMatrix$, or $
\Lambdabsf \ybsf^{\star}_i$ and $\ybsf^{\star}_i
$ are parallel. 
We can rule out the former: because $\lambda_j \neq 0$ for all $j \in \Ic$, then $\Lambdabsf \ybsf^{\star}_i = \ZeroMatrix$ implies that $\langle \ybsfs_i, \ybsfs_i \rangle= 0$, 
which contradicts the fact that $\pr{\ybsfs_1, \ldots, \ybsfs_n}\in \Oprime{n}$ and therefore $\langle \ybsfs_i, \ybsfs_i \rangle= 1$. 
Thus,
\begin{equation}
    \Lambdabsf \ybsf^{\star}_i = \alpha \ybsf^{\star}_i,
\label{eq:ystarform}
\end{equation}
for some $\alpha \in \Rbb$.

\item Let $\pr{\ebsf_1,\ldots, \ebsf_n}$ be the canonical set of basis vectors of $\Rbb^n$. 
That is, for each $i\in \Ic$, $\ebsf_i =\pr{\delta_{ik}}_{k\in \mathcal{I}}$. 
Since $\lambda_i \neq \lambda_j$ for $i\neq j$ and $\langle \ybsf_i^{\star}, \ybsf_i^{\star}\rangle =1$, then \eqref{eq:ystarform} implies that for each $i\in \Ic$ the vector $\ybsf^{\star}_i$ has to be equal to $\pm \ebsf_j$ where $j \in \Ic$. 
The constraint that $\langle \ybsf_i^{\star}, \ybsf_j^{\star}\rangle =0$ for $i\neq j$ implies that two different $\ybsf^{\star}_i$ cannot correspond to the same canonical basis vector. 
Hence, it follows that $\pr{\ybsfs_1,\ldots, \ybsfs_n}$ has to be of the form $\pr{\pm \ebsf_{\sigma\ag{1}},\pm \ebsf_{\sigma\ag{2}},\ldots, \pm \ebsf_{\sigma\ag{n}}}$, where $\sigma\ag{\cdot}$ is a permutation on $\Ic$. 

\item Thus, the stationary points of $\tilde{g}\ag{\cdot}$, which includes the global maximizer $\pr{\ybsfs_1,\ldots, \ybsfs_n}$, belong to the finite set
\begin{equation}
\left\{\pr{\pm \ebsf_{\sigma\ag{1}},\pm \ebsf_{\sigma\ag{2}},\ldots, \pm \ebsf_{\sigma\ag{n}}}\SUCHTHAT\text{ $\sigma\ag{\cdot}$ is a permutation on $\Ic$}\right\}.
\label{eq:Ysol}
\end{equation}
The image of this set under $\tilde{g}\ag{\cdot}$ is
\begin{equation}
 \left\{ \sum_{i=1}^n \mu_i \lambda_{\sigma\ag{i}}\SUCHTHAT \text{ $\sigma\ag{\cdot}$ is a permutation on $\Ic$}
\right\},
\label{eq:gatSta}
\end{equation}
which must contain the maximum value of $\tilde{g}\ag{\cdot}$. The maximal element of this finite set is $\sum_{i}\lambda_i\mu_i$\footnote{This is because $\pr{\lambda_i}_{i \in \Ic}$ and $\pr{\mu_i}_{i \in \Ic}$ are descending lists of numbers. To elaborate further note  that if $a_1 \geq a_2$ and $b_1 \geq b_2$, then $(a_1 - a_2)(b_1 - b_2) \geq 0$ implies that $a_1 b_1 + a_2 b_2 \geq a_1 b_2 + a_2 b_1$.}. 

\item The maximal element of the set \eqref{eq:gatSta} corresponds to the identity permutation. 
Hence, it follows from \eqref{eq:Ysol} that the set of all $\Ytbsf^{\star}$, i.e., the set of all global maximizers of $\tilde{g}\ag{\cdot}$, equals
\begin{equation}
\tilde{S}^{\star}_n=\left\{\pr{\pm \ebsf_{1},\pm \ebsf_{2},\ldots, \pm \ebsf_{n}}\right\},
\label{eq:Ysol2}
\end{equation}
among which is $\widetilde{\IdentityMatrix}_n$.
\end{enumerate}

\item\label{thm:lambda-m-n:step2} In this step we relax the conditions we required on $\lambda_i$ in Step~\ref{thm:lambda-m-n:step1}, 
and assume only that $\lambda_i$ are non-increasing.
However, we still require $\mu_i$ to be strictly decreasing and non-zero. 
\begin{enumerate}
\item For any $\varepsilon > 0$, define $\lambda_i(\varepsilon) = \lambda_i + \varepsilon (n + 1 - i)$, and let $\Lambdabsf_{\varepsilon} = \DIAG_n\ag{\lambda_1(\varepsilon),\lambda_2(\varepsilon),\ldots}$. 
When $\varepsilon$ is sufficiently small, we guarantee that $\lambda_i(\varepsilon) > \lambda_{i+1}(\varepsilon)$ and $\lambda_i(\varepsilon) \neq 0$ hold for all $i$. 
That is, for $\varepsilon$ sufficiently small $\lambda_i\pr{\varepsilon}$ satisfy the same conditions that $\lambda_i$ did in Step~\ref{thm:lambda-m-n:step1}. 
\item Define $\tilde{g}_\varepsilon\ag{\cdot}:\Oprime{n}\to \Rbb$ as
\begin{equation}
   \tilde{g}_\varepsilon\ag{\pr{\ybsf_1,\ldots, \ybsf_n}}= \sum_{i=1}^n \mu_i \langle \ybsf_i,  \Lambdabsf_{\varepsilon} \ybsf_i\rangle.
   \label{eq:fe}
\end{equation} 
It follows from \eqref{eq:f} that
 $$
\begin{aligned}
\tilde{g}_\varepsilon\ag{\Ytbsf}
&= \tilde{g}\ag{\Ytbsf} + \varepsilon \sum_{i=1}^n \mu_i \langle \ybsf_i, \DIAG_n\ag{n,n-1,\ldots,1} \ybsf_i\rangle,
\end{aligned}
$$
where we recall that $\Ytbsf := \pr{\ybsf_1,\ldots, \ybsf_n}$. 
For the second term, we have
$$
\varepsilon\sum_{i=1}^n \mu_i \langle \ybsf_i, \DIAG_n\ag{n,n-1,\ldots,1} \ybsf_i\rangle \leq \varepsilon \sum_{i=1}^n | \mu_i | n \langle \ybsf_i, \ybsf_i \rangle= \varepsilon n \sum_{i=1}^n | \mu_i |,
$$
which is a constant that is independent of $\Ytbsf$. 
Hence,
\begin{equation}
\left| \tilde{g}_\varepsilon\ag{\Ytbsf} - \tilde{g}\ag{\Ytbsf} \right| \leq \varepsilon n \sum_{i=1}^n | \mu_i |.
\label{ineq:geg}
\end{equation}
\item\label{thm:lambda-m-n:step2:c}
Recall that $\Itbsf_n$ belongs to the set $\tilde{S}^{\star}_n$ defined in \eqref{eq:Ysol2}. 
Then, for any $\Ytbsf$
\begin{align}
\tilde{g}\ag{\,\Itbsf_n} - \tilde{g}\ag{\Ytbsf} = \pr{\tilde{g}\ag{\,\Itbsf_n} - \tilde{g}_\varepsilon\ag{\,\Itbsf_n}} + \pr{\tilde{g}_\varepsilon\ag{\,\Itbsf_n} - \tilde{g}_\varepsilon\ag{\Ytbsf}} + \pr{\tilde{g}_\varepsilon\ag{\Ytbsf} - \tilde{g}\ag{\Ytbsf}}.
\label{eq:addsub}
\end{align}
Using \eqref{ineq:geg} for the first and third terms, we bound the previous equation from below:
\begin{align}
\tilde{g}\ag{\,\Itbsf_n} - \tilde{g}\ag{\Ytbsf} &\geq \pr{\tilde{g}_{\varepsilon}\ag{\,\Itbsf_n} - \tilde{g}_\varepsilon\ag{\Ytbsf}} - 2 \varepsilon n \sum_{i=1}^n | \mu_i |.
\label{ineq:addsub1}
\end{align}
Note that  $\tilde{g}_{\varepsilon}\ag{\cdot}$ satisfies all the conditions imposed on $g\ag{\cdot}$ in Step~\ref{thm:lambda-m-n:step1}, and hence $\tilde{g}_{\varepsilon}\ag{\cdot}$ achieves its global maxima at $\Itbsf_n$. 
That is,
$\tilde{g}_{\varepsilon}\ag{\,\Itbsf_n} \ge \tilde{g}_\varepsilon\ag{\Ytbsf} $ for all $\Ytbsf\in \Oprime{n}$. 
Then it follows from \eqref{ineq:addsub1} that
\begin{align}
\tilde{g}\ag{\,\Itbsf_n} - \tilde{g}\ag{\Ytbsf}
&\geq - 2 \varepsilon n \sum_{i=1}^n | \mu_i |.
\end{align} 
Since $\varepsilon > 0$  can be arbitrarily small, the above inequality implies that
\begin{equation}
\tilde{g}\ag{\,\Itbsf_n} - \tilde{g}\ag{\Ytbsf}\geq 0.
\label{ineq:bound}
\end{equation}
We could have chosen any $\Ytbsf^\star \in \tilde{S}^\star_n$ and arrived at the same conclusion; thus we can also conclude that $\tilde{S}^{\star}_n$ is the set of global maximizers of $\tilde{g}\ag{\cdot}$ even in this more general case. 
\end{enumerate}
\newcommand{\Dbsfe}{\Dbsf\pr{\epsilon}}
\item Finally, we relax the conditions on $\mu_i$ as well. 
That is, we consider a $\tilde{g}\ag{\cdot}$ in which  we only require that $\lambda_i$ and $\mu_i$  each be a non-increasing sequence. 
The procedure we follow is similar to the one in the previous step. 
\begin{enumerate}
    \item As we did for $\lambda_i$ in Step~\ref{thm:lambda-m-n:step2}, we construct a perturbed sequence $\mu_i(\varepsilon)$, where $\mu_i(\varepsilon):=\mu_i+\varepsilon(n+1-i)$. 
    Again, for $\varepsilon$ sufficiently small  $\mu_i(\varepsilon)$ will be a strictly decreasing sequence and $\mu_i(\varepsilon)\neq0$ for each $i\in \Ic$. 
    We define $\Dbsf(\varepsilon)=\DIAG_n\ag{\mu_1(\varepsilon), \ldots,\mu_n(\varepsilon)}$. 
    \item Let $\tilde{g}_\epsilon\ag{\cdot}:\Oprime{n}\to \Rbb$,
\begin{equation}
   \tilde{g}_\varepsilon\ag{\pr{\ybsf_1,\ldots, \ybsf_n}} = \fprod{ \Lambdabsf, \Ybsf \Dbsf(\varepsilon) \Ybsf^{\Tsf}},
   \label{eq:fe1}
\end{equation}
where $\Ybsf=\texttt{iso}^{-1}\ag{\ybsf_1,\ldots, \ybsf_n}$. 
Similar to the previous step, it can be shown with \eqref{eq:f} that
\begin{align}
\tilde{g}_\epsilon\ag{\Ytbsf} = \tilde{g}\ag{\Ytbsf}
+ \sum_{i=1}^n \varepsilon (n + 1 - i) \langle \ybsf_i, \Lambdabsf \ybsf_i \rangle.
\label{eq:geexpanded}
\end{align} 
For the second term,
\begin{equation}
\left|\langle \ybsf_i, \Lambdabsf \ybsf_i \rangle \right| \leq \max_j |\lambda_j| \langle \ybsf_i, \ybsf_i \rangle = \max_j |\lambda_j|.
\end{equation}
Thus,
\begin{align}
\left| \tilde{g}_\epsilon\ag{\Ytbsf} - \tilde{g}\ag{\Ytbsf} \right| &\leq \varepsilon \max_j |\lambda_j| \sum_{i=1}^n (n +1 - i), \notag \\
&= \frac{1}{2} n (n + 1) \varepsilon \max_j |\lambda_j|.\label{ineq:geg2}
\end{align}

\item Using calculations similar to those in Step~(\ref{thm:lambda-m-n:step2}.c) we get from \eqref{ineq:geg2}
\begin{align}
\tilde{g}\ag{\,\Itbsf_n} - \tilde{g}\ag{\Ytbsf}&\geq \pr{\tilde{g}_{\epsilon}\ag{\,\Itbsf_n} - \tilde{g}_\epsilon\ag{\Ytbsf}} - n (n + 1) \varepsilon \max_j |\lambda_j|,
\label{ineq:5addsub1}
\end{align}
where $\Itbsf_n$ is an element of $\tilde{S}^{\star}_n$. 
Note that  $\tilde{g}_{\epsilon}\ag{\cdot}$ satisfies all the conditions imposed on $\tilde{g}\ag{\cdot}$ in Step~\ref{thm:lambda-m-n:step2}, hence $\tilde{g}_{\epsilon}\ag{\cdot}$ achieves its global maximum at $\Itbsf_n$. 
That is,
$\tilde{g}_{\epsilon}\ag{\,\Itbsf_n} \ge \tilde{g}_\epsilon\ag{\Ytbsf} $ for all $\Ytbsf\in \Oprime{n}$. 
Thus,
\begin{align}
\tilde{g}\ag{\,\Itbsf_n} - \tilde{g}\ag{\Ytbsf}&\geq  - 2 \epsilon n \sqrt{n} \sum_{i=1}^n | \lambda_i |.
\label{ineq:5addsub2}
\end{align} 
Using arguments similar to those we used towards the end of Step~(\ref{thm:lambda-m-n:step2}.c) it can be concluded from \eqref{ineq:5addsub2} that even this more general $\tilde{g}\ag{\cdot}$ achieves its maximum value on the set $\tilde{S}^{\star}_n$.
\end{enumerate}

\end{enumerate}
\end{proof}

\bibliography{references.bib}

\end{document}